\documentstyle[amstex,amscd,12pt]{article}
\makeatletter \oddsidemargin5pt \evensidemargin5pt \textwidth 16cm
\date{}

\newtheorem{theorem}{Theorem}[section]

\newtheorem{proposition}[theorem]{Proposition}

\newtheorem{lemma}[theorem]{Lemma}

\newcommand{\z}{{\Bbb Z}}

\newcommand{\N}{{\Bbb N}}

\newcommand{\st}{{\rm st}}
\newcommand{\Int}{{\rm Int}}
\newcommand{\cl}{{\rm cl}}

\newcommand{\no}{$\rm  {N\ddot{o}beling}$ }
\newcommand{\noo}{$ {N\ddot{o}beling}$ }

\newcommand{\lo}{\longrightarrow}

\begin{document}

\title{Characterizing  N\"{o}beling spaces\\
}

\author{Michael  Levin}

\maketitle
\begin{abstract}
It  is shown that \no spaces are uniquely determined by the
universal  extension and embedding properties.
\bigskip
\\
{\bf Math. Subj. Class.:} 55M10, 54F45.
\end{abstract}
\begin{section}{Introduction}\label{sec1}

      All spaces are assumed to be separable metrizable.
 A complete $n$-dimensional metric space $X$ is  said to be an
 $n$-dimensional \no space if
 the following properties are satisfied:

 (i) $X $ is an absolute extensor in dimension $n$,
  that is every map
 $f: A \lo X$ from a closed subset $A$ of a space $Y$ of
 $\dim \leq n$
  extends over $Y$;

 (ii) every map $f : Y \lo X$ from a complete metric space
 $Y$ of $\dim \leq n$ can be arbitrarily closely
 approximated by a closed embedding, that is  for every open
 cover  $\cal U$ of $X$ there is a closed embedding $g : Y \lo X$
 which is $\cal U$-close to $f$
 ($\cal U$-close means that  for every $y \in Y$ there is
 an element of $\cal U$ that contains  both $f(y)$ and $g(y)$).
 \\

  The  characterization theorem for \no spaces says that:
   \begin{theorem}\label{char}
   Every two \noo spaces of the same dimension are
   homeomorphic.\footnote{ This theorem was also obtained by
      A. Nagorko  \cite{nagor} using a different approach.
       The results of \cite{nagor} were announced in
       July 2005  in Bedlewo, Poland.}
   \end{theorem}

   This result is well-known in dimension $0$. One-dimensional
   \no spaces were characterized in \cite{klt}. Our main goal  is
    to prove the characterization theorem in all
   dimensions. This paper is a continuation of \cite{lev} where
   key ingredients  and techniques for proving Theorem \ref{char}
   were established. The general pattern of the proof
   of Theorem \ref{char} follows Toru\'nczyk's and Bestvina's
   proofs of the  characterization theorems for
   the Hilbert cube \cite{tur1}, the Hilbert space \cite{tur2}
   and the Menger universal
   compacta  \cite{be}. Let us
   remind some  facts from \cite{lev} and
    outline key steps
   of the proof of Theorem  \ref{char}.

 Throughout the paper   a manifold means a manifold with
 (possibly empty) boundary and
   a triangulated space   means a locally finite simplicial
  complex which we identify with the underlying space.
  For a triangulated space we consider only triangulations
  compatible with the PL-structure of the space.
  All triangulated manifolds are assumed to be combinatorial.
  We will  work   with the following model of   \no spaces.

  By a rational map between triangulated spaces we mean a PL-map
  that sends points with rational barycentric coordinates
  to points with rational barycentric coordinates.
   Two triangulations  of a space   are said to be rationally
  equivalent if the identity map is a rational map with
  respect to these triangulations (it is easy to see that
  if a PL-homeomorphism is rational in one direction then it is
  rational in the opposite direction as well).
  Let $M$ be a triangulated space. Every triangulation of $M$
  which is rationally equivalent to the given triangulation
  of $M$ is said to be a rational triangulation and
  the class of all rational triangulations is said to be
  the rational structure of $M$.
   Denote by
  $M(k)$ the subspace of $M$ which is the complement of
  the union of all the triangulated spaces of $\dim \leq k$
  which are rationally  embedded in $M$.

  Let us state the following important fact  leaving its proof to the
  reader.

 \begin{theorem}
 \label{model} Let $M$ be a triangulated $m$-dimensional manifold,
  let $k\geq 0$ be an integer and
 let $n=m-k-1$.  If $M$ is $(n-1)$-connected and
 $m\geq 2n+1$ then $M(k)$ is an $n$-dimensional \noo space.

 \end{theorem}

  A space  $M(k)$ satisfying  the assumptions of Theorem \ref{model}
 is called a  \no space modeled on a triangulated manifold.

 A  subset $A$  of a space  $X$ is called a $Z$-set if
 $A$ is closed in $X$ and
 the identity map of $X$ can be arbitrarily closely approximated
 by a map $f : X \lo X$ with $f(X) \cap A =\emptyset$.
 Note that if $X$ is an $n$-dimensional \no space modeled
 on a manifold $M$
  and $A \subset X$ is a $Z$-set in $X$ then
 $X \setminus A$ is also an $n$-dimensional
  \no space modeled on the manifold $N=M \setminus$ the closure
  of $A$ in $M$ (the rational structure of $N$ is defined such that
  the inclusion is a rational map).
 The following version of a $Z$-set unknotting theorem was proved in \cite{lev}.

 \begin{theorem}
 \label{unknot}{\rm \bf(Unknotting Theorem \cite{lev})}
 Let $X_1$ and $X_2$ be $n$-dimensional \noo spaces and let $A_1$
 and $A_2$ be $Z$-sets in $X_1$ and $X_2$ respectively
 such that $X_1\setminus A_1$ and $X_2 \setminus A_2$  are
 homeomorphic to  $n$-dimensional \noo spaces modeled
 on triangulated manifolds.
 If $A_1$ and $A_2$ are homeomorphic
 then any homeomorphism  between $ A_1$ and $ A_2$ can be extended to
 a homeomorphism between $X_1$ and $X_2$.
 \end{theorem}

 A subset  of a space $X$ is called a $\sigma$-$Z$-set if it is
 a countable union of $Z$-sets. In Section \ref{sec-resol} we prove
 the following version of a  resolution theorem which establishes
  a connection between general \no spaces
 and \no spaces modeled on triangulated manifolds.

 \begin{theorem}
 \label{resol}{\rm \bf{(Resolution Theorem})}
 Let $Y$ be an $n$-dimensional  \noo space. Then there is
 a subspace $X \subset Y$ such that $X$ is homemorphic
 to an $n$-dimensional \noo space modeled on a triangulated
 manifold and $Y\setminus X$ is a $\sigma$-$Z$-set in $Y$.
 \end{theorem}

  In Section \ref{sec-shrink} we prove the following version of a shrinking
  theorem which is the most important step in proving Theorem \ref{char}.

 \begin{theorem}
 \label{shrink}
 {\rm \bf{(Shrinking Theorem})}
 Let $X$ be a subspace of an $n$-dimensional  \noo space $Y$ such that
 $X$ is homeomorphic to an $n$-dimensional \noo
 space modeled on a triangulated manifold
 and $Y \setminus X$ is a $\sigma$-$ Z$-set in $Y$.
  Then for every
 covers ${\cal V}_X$ and ${\cal V}_Y$ of $X$ and $Y$
 by  sets open in $X$ and $Y$ respectively
 there is a homeomorphism $h : X \lo X$ such that
 for every $x \in X$ the point  $h(x)$
 is ${\cal V}_Y$-close to $x$
    and for every
   $y \in Y$ there are a neighborhood $G $ of $y$ in $Y$
   and a set $V \in {\cal V}_X$ such that $h^{-1}(G\cap X) \subset V$.
   \end{theorem}

 Theorems \ref{unknot}, \ref{resol} and \ref{shrink} imply
 Theorem \ref{char}. Indeed,
   assume that $X$ and $Y$ satisfy the assumptions of Theorem
  \ref{shrink}.
  Since $Y\setminus X$ is a $\sigma$-$Z$-set in $Y$ we have that
  $X$ is dense in $Y$. Then
   by Bing's Shrinking Criterion  and Theorem \ref{shrink}
   the inclusion of $ X$ into $ Y$ can be arbitrarily closely
   approximated by a homeomorphism between $X$ and $Y$.
   Thus by Theorem \ref{resol} we get that every $n$-dimensional
    \no space is homeomorphic
   to an $n$-dimensional \no space modeled on a triangulated
   manifold. By Theorem \ref{unknot} every two
   $n$-dimensional  \no spaces modeled on triangulated
   monifolds are homeomorphic, and  Theorem \ref{char} follows.

    As we already mentioned  this paper is a continuation of
     \cite{lev}. Our presentation heavily relies on \cite{lev}
     and  we assume that the reader is closely familiar
     with the results and the constructions of \cite{lev}.

    \end{section}

     \begin{section}{Unknotting and resolution theorems}\label{sec-resol}

       For  an open subset $V$ of a dense subset $D$
     of a space $X$,
     by the extension of $V$ to $X$
      we understand
       the largest open set $V'$ in $X$
     such that $V=D\cap V'$. Similarly,   by  the extension to $X$ of
      a collection     ${\cal V}$
     of open subsets of $D$ we mean the collection
     consisting of the extensions to $X$  of the elements of $\cal V$.

     For a collection $\cal V$ of subsets  of $X$  denote
     $\st^0({\cal V})={\cal V}$,
     $\st({\cal V})=\st^1({\cal V})=\st({\cal V}, {\cal V})$ and
     define
     by induction $\st^{n+1}({\cal V})=\st( \st^n({\cal V}))$.
     We write ${\cal  V} \prec {\cal U}$ if
     ${\cal V}$
      refines a collection $\cal U$ of subsets of $X$.


     A collection    $\cal V$  of subsets of $X$ is said to be
     an $(n-1)$-refinement
     of  a collection $\cal U$ of subsets of $X$,
     written ${\cal V} \prec_{n-1} \cal U$, if
      for every $V \in \cal V$ there is $U \in \cal U$ such that
      $V \subset U$ and the inclusion of $V$ into $U$ induces the
      zero-homomorphism of the homotopy groups in $\dim \leq n-1$

     A map $f : X \lo Y$ is said to be $UV^{n-1}$ if $f(X)$ is dense
     in $Y$ and for every $y \in Y$ and every neighborhood $U$
     of $y$ there is a smaller neighborhood $V$ of
     $y$
       such that the inclusion of
     $f^{-1}({ V})$ into
      $f^{-1}({ U})$ induces the zero-homomorphism of
      the homotopy groups in $\dim \leq n-1$.

     A map $f : X \lo Y$ is said to be a $Z$-embedding if
     $f$ is a closed embedding of $X$ into $Y$ and
     $f(X)$ is a $Z$-set in $Y$.

     The following proposition whose proof is left to the reader
     presents one of the basic properties of $UV^{n-1}$-maps
     on \no spaces.

     \begin{proposition}
     \label{p2.1}
     Suppose that $f : X \lo Y$ is a $UV^{n-1}$-map from
     an $n$-dimensional  \noo space
     $X$  and  ${\cal U}$ is an open cover of $Y$.
     Then there is an open cover ${\cal V}$ of $Y$ having
     the property:
     given a closed subset $A$ of a complete space $B$ of
     $\dim     \leq n$,
     a $Z$-embedding  $\phi_A : A \lo X$ of $A$ into $X$
     and a map $\psi_B : B \lo Y$
     such that $\psi_B|A$ and $f \circ \phi_A$ are $\cal V$-close,
        the map
     $\phi_A$ extends to a $Z$-embedding $\phi_B : B \lo X$
     of $B$ into $X$  such that
     $\psi_B$ and $f\circ \phi_B$ are ${\cal U}$-close.
     \end{proposition}

   The following version of a $Z$-set unknotting
    theorem   generalizes Theorem \ref{unknot}.

    \begin{theorem}
    \label{unknot1}
    Let $f_i : X_i \lo Y, i=1,2$ be $UV^{n-1}$-maps from
    $n$-dimensional \noo spaces $X_i$ modeled on triangulated
    manifolds  and
    let $\cal U$
    be an open cover of $Y$. Then there is an open cover
    $\cal V$ of $Y$ such that for every pair of homeomorphic
     $Z$-sets    $A_i \subset X_i$ and a homeomorphism
     $h_A : A_1 \lo A_2$ such that the maps $f_1|A_1$ and
     $f_2\circ h_A$ are ${\cal V}$-close, the homeomorphism
     $h_A$ extends to a homeomorphism $h : X_1 \lo X_2$
     such that the maps $f_1$ and
     $f_2\circ h$ are ${\cal U}$-close.
     \end{theorem}
    For proving Theorem \ref{unknot1} we need the following
    simple modification of Proposition 3.1 of \cite{lev}.
       \begin{proposition}
  \label{p2.2}
  Let $X$ be an $n$-dimensional \noo space,
    $A$  a $Z$-set in $X$,  $\cal C$
    a cover of $X \setminus A$ that properly approaches $A$
   and  $\cal W$  an open cover of $X$ such that
  $\cal C$ is an $(n-1)$-refinement of $\cal W$.
   Then for every $C \in \cal C$  there is an open set
   $C\subset V_C \subset X \setminus A$ such that
   the inclusion $C \subset V_C$ induces the zero-homomorphism
   of the homotopy groups in $\dim \leq n-1$,
    ${\cal V} =\{ V_C : C \in {\cal C} \}$ properly approaches
    $A$
     and $\cal V$
    refines
    $\st^3{\cal W}$.
   \end{proposition}
   {\bf  Proof.} The only change in the proof of
   Proposition 3.1 of \cite{lev} that we need to make
   is to assume that  the sets $G$ and $G_C$ are contained
   in elements  of $\cal W$, the maps $e_n$ are $\cal W$-close
   to the identity map of $X$ and $\cal U$ refines $\cal W$.
   Then $\{ Y_C : C \in {\cal C}\}$ refines $\st^3{\cal W}$
   and therefore $\cal V$ can be chosen so that
   $\cal V$
    refines  $\st^3{\cal W}$.
        \hfill $\Box$
 \\\\
     {\bf  Proof of Theorem \ref{unknot1}.}
     The proof is similar to the proof of the unknotting
     theorem in \cite{lev}.
     The homeomorphism between $X_1$ and $X_2$ constructed
     in the proof of Theorem 1.2 of \cite{lev} will have the
     required properties if  the partitions ${\cal P}_1$
     and ${\cal P}_2$ are constructed so that
     for every $P \in {\cal P}_1$ there is an element
     of $\cal U$ which contains both  $f_1(P \cap X_1)$ and
     $f_2(\mu(P)\cap X_2)$. In order to gain this additional
     control on ${\cal P}_1$ and ${\cal P}_2$ we need
     to make the following
     adjustments in
     the proof  of
     Proposition 3.2 of \cite{lev}.

     To avoid any possible
     confusion with the notation
     of the proof of
     Proposition 3.2 of \cite{lev} let us re-denote
     $\cal V$ and $\cal U$ by ${\cal V}_Y$ and
     ${\cal U}_Y$, $f_1$ and $f_2$ by
     $f_{X_1}$ and $f_{X_2}$ and $h$ by $h_X$  respectively.
      Now we can  adopt the notation of
         Proposition 3.2 of \cite{lev}.

   Let $\omega$ be the number of times that we use the
     constructions 2.7  and 2.9 of \cite{lev}  for improving
     connectivity of ${\cal P}_2$ in the proof of 3.2
     of \cite{lev}. Note that $\omega$ is finite and
     depends only on the dimension  $n$.
     Because $f_{X_i}$ is $UV^{n-1}$ we can choose
      a finite
     sequence ${\cal Y}^j,  j=1, \dots,2\omega$ of
     open covers of $Y$ such that ${\cal Y}^{2\omega}$
     refines ${\cal U}_Y$ and for
     ${\cal X}^j_i = f_{X_i}^{-1}({\cal Y}^j)$
     we have that $\st^5({\cal X}^j_i)$
     is an $(n-1)$-refinement of
      ${\cal X}^{j+1}_i$,
      $j=1, \dots,2\omega-1$. Set ${\cal V}_Y ={\cal Y}^1$.

     Fix a pair of $Z$-sets $A_i \subset X_i$
     and a homeomorphism $h_A : A_1 \lo A_2$
     such that $f_{X_2}\circ h_A$ and  $f_{X_1}|A_1$ are
      ${\cal Y}^1$-close.
      Recall that
     in  Proposition 3.2 of \cite{lev}
       $X_i$ is considered as a subspace of a space
     $Y_i$ that can be represented as $Y_i=M_i \cup X_i$
     such that $M_i$ is a triangulated manifold,
     $X_i \setminus A_i = M_i(k_i)$, $A_i = Y_i \setminus M_i$
     and $A_i$ is closed in $Y_i$.  Following
     the proof of Proposition 3.2 of \cite{lev} we
      make
     the following adjustments.

      We can assume that $M'_1$ and a partition ${\cal P}_1$
      of $M'_1$ are chosen so that
      the extension ${\cal Y}^j_1$ of ${\cal X}^j_1$ to $Y_1$
      covers $M'_1$
      and
      ${\cal P}_1$
      refines ${\cal Y}^1_1$. Since $f_{X_i}$ is
      $UV^{n-1}$, we can  assume
      by Proposition      \ref{p2.1} that the map
      $g : X_2 \lo X_1$ is chosen so that $f_{X_2}: X_2 \lo Y$ and
      $f_{X_1}\circ g : X_2 \lo Y$ are ${\cal Y}^1$-close.
      Then, it is easy to see that $M'_2$ and  the initial
      partition  ${\cal P}_2$ of $M'_2$ and  can be constructed so that
      the extension ${\cal Y}^j_2$
      of  ${\cal X}^j_2$ to $Y_2$
       covers  $M'_2$ and the one-to-one correspondence
       $\mu : {\cal P}_1 \lo {\cal P}_2$ is ${\cal Y}^1$-agreeable with
      $f_{X_1}$ and $f_{X_1}$.
      The last property means that for every $P\in {\cal P}_1$
      there is an element of ${\cal Y}^1$ which contains
      both $f_{X_1}(P\cap X_1 )$  and $f_{X_2}(\mu(P) \cap X_2)$.
      Note that since every map from a sphere of $\dim \leq n-1$
      into an open subset $G$ of $Y_i$ can be homotoped
      inside $G$ into $G\cap X_i$ we have that
      $\st^5({\cal Y}^j_i)$
     is an $(n-1)$-refinement of
      ${\cal Y}^{j+1}_i$,
      $j=1, \dots,2\omega-1$.

      Our next step is to analyze the procedure
      of  improving   connectivity of
      ${\cal      P}_2$.
      Recall that, in order to simplify the notation, in the beginning
      of each step of the procedure we replace  $M_2$ and
      ${\cal P}_2$ by  the output $M'_2$ and ${\cal P}'_2$ of
      the previous step (the modifications of
      $M_2$ and ${\cal P}_2$) and  we also use  $\mu$ to
      denote the one-to-one correspondence
      $\mu : {\cal P}_1 \lo {\cal P}'_2$ which is the  composition
      of $\mu : {\cal P}_1 \lo {\cal P}_2$
      with the natural correspondence between
      ${\cal P}_2$ and its modification ${\cal P}'_2$.

      Thus we replace  $M_2$ by $M'_2$ for our initial partition
      ${\cal P}_2$ and now we
      assume that at some step of the procedure
      $M_2$ and ${\cal P}_2$ are already constructed so that
      $\st({\cal P}_2)$ refines ${\cal Y}^j_2$ and
      the one-to-one correspondence
      $\mu : {\cal P}_1 \lo {\cal P}_2$ is ${\cal Y}^j$-agreeable
      with $f_{X_2}$ and $f_{X_1}$.

       The following adjustments  should be made
      in the part
      the proof of
     Proposition 3.2 of \cite{lev} where
     the construction 2.7 of \cite{lev} is applied.
     Choose $\cal C$ so that $\cal C$ refines ${\cal Y}_2^j$.
      Then by Proposition \ref{p2.2} the cover $\cal V$ of
      $M_2(k_2)$ can be chosen so that ${\cal V}$ refines
      $\st^3({\cal X}_2^{j+1})$. It implies that
      ${\cal V}$ refines ${\cal  Y}^{j+1}_2$ and hence
      ${\cal W}=\st(C, {\cal V})$ refines $\st({\cal Y}^{j+1}_2)$.
      Then the needed  modification of ${\cal P}'_2$ which is
      (the output of 2.7 of
      \cite{lev}) can be constructed so that it will refine
      $\st^2({\cal W})$ and hence it will refine
      $\st^3({\cal     Y}_2^{j+1})$ and, as a result, will
      refine  ${\cal Y}^{j+2}_2$ as well.

      Now we will adjust
      the part of
      the proof of
     Proposition 3.2 of \cite{lev} where
     the construction 2.9 of \cite{lev} is applied.
          Once again we assume that $M_2$ and ${\cal P}_2$
           are constructed so that
      $\st({\cal P}_2)$ refines ${\cal Y}^j_2$ and
      the one-to-one correspondence
      $\mu : {\cal P}_1 \lo {\cal P}_2$ is ${\cal Y}^j$-agreeable
      with $f_{X_2}$ and $f_{X_1}$.
      Then ${\mu^{-1}(\st(\cal P}_2))$ refines $\st({\cal Y}^{j}_1)$
      and therefore by Proposition \ref{p2.2} the cover
      ${\cal V}$ of $M_1(k_1)$ can be chosen so that
      ${\cal V}$ refines ${\cal Y}^{j+1}_1$. Then
      ${\cal H}=\st({\cal V}, \mu^{-1}(\st({\cal P}_2)))$ refines
      $\st({\cal Y}^{j+1}_1)$ and hence $\mu({\cal H})$
      refines ${\cal Y}^{j+1}_2$. Now  ${\cal W}$ can be chosen so
      that ${\cal W}$ refines ${\cal Y}^{j+1}_2$ and therefore
      the modification of ${\cal P}_2$ will refine
      $\st^2({\cal  Y}^{j+1}_2) $ and, as a result, will refine
       ${\cal Y}^{j+2}_2$ as well.

        In both   2.7 and
         2.9 of \cite{lev} the modification ${\cal P}'_2$
         of ${\cal P}_2$ can be constructed such that
        for every element $P\in{\cal P}_2$, its
      modification $P'\in{\cal P}'_2$
      intersects $P$. Then  we get
      that $\mu : {\cal P}_1 \lo {\cal P}'_2  $
      is ${\cal Y}^{j+2}$-agreeable with $f_{X_1}$
      and $f_{X_2}$ because $\mu$ was ${\cal Y}^j$-agreeable
      before the constructions. In addition,
      applying a homeomorphism $h : M_2\cup A_2  \lo M_2 \cup A_2$
      sufficiently close to the identity map
      as it is described in the proof of Proposition 3.2
      of      \cite{lev}
       we can replace
      $M'_2$ and  ${\cal P}'_2$ by  $h(M'_2)$ and  $h({\cal P}'_2)$
      respectively and  assume that ${\cal P}'_2$
      is a rational decomposition of $M'_2$.

      Thus denoting the final modification of ${\cal P}_2$ again
      by
      ${\cal P}_2$ we get that $\mu : {\cal P}_1 \lo {\cal P}_2$
      is ${\cal Y}^{2\omega}$-agreeable
      with $f_{X_1}$
      and $f_{X_2}$.  Recall that ${\cal Y}^{2\omega}$ refines
      ${\cal U}_Y$. Now we can construct a homeomorphism
      $h_X: X_1 \lo X_2$ which carries $P \cap X_1$
      onto $\mu(P)\cap X_2$ for every $P \in {\cal P}_1$.
      Then $f_{X_2}\circ h_X$ and $f_{X_1}$ are ${\cal U}_Y$-close
      and the theorem is proved.
            \hfill $\Box$

     \begin{proposition}\label{p2.3}
     Let $f : X \lo Y$ be a $UV^{n-1}$-map from
     an $n$-dimensional \noo space $X$ modeled on a triangulated
     manifold to a complete space $Y$
     and let $g_A : A \lo Y$ be a map from a $Z$-set
     $A \subset X$. Then there is a $UV^{n-1}$-map
     $g : X \lo Y$ such that $g|A=g_A$.
     \end{proposition}
      {\bf  Proof.}
       Fix a complete  metric in $Y$
      with distances bounded by $1/8$. Set
      ${\cal Y}_0={\cal Y}_1={\cal Y}_2={\cal Y}_3=$the
      trivial cover $Y$ consisting of only one set $Y$.
      Using
      Theorem \ref{unknot1} choose  a sequence
      of open covers ${\cal Y}_i$ of $Y$ such that

      (1) mesh${\cal Y}_i\leq 1/2^i$;

      (2) $\st^2{\cal Y}_{i+1}\prec {\cal Y}_i$;

     (3)  $f^{-1}({\cal Y}_{i+1})
       \prec_{n-1}f^{-1}({\cal Y}_{i})$;

    (4)  for every homeomorphism $h_B: B_1 \lo B_2$ of
    $Z$-sets in $X$ such that $f|{B_1}$ and $f\circ h_B$
    are ${\cal Y}_{i+1}$-close, $h_B$ extends to
    a homeomorphism $h : X \lo X$ such that
    $f$ and $f\circ h$ are ${\cal Y}_i$-close.

    Set $A_0=A$ and $h_0=$the identity map of $X$.
     We are going to construct for every $i$
    a homeomorphism $h_i: X \lo X$ such that

    (5) $f\circ h_i$ and $f$ are ${\cal Y}_i$-close;

    (6) $f \circ h^{i}|A$ and $g_A$ are ${\cal Y}_{i+3}$-close
      where $h^{i}= h_{i} \circ\dots \circ h_0 : X \lo X$.
     \\    Assume that the construction is completed up to the index $i$.
    Proceed to $i+1$ as follows. By
    Proposition \ref{p2.1} choose
     a $Z$-embedding $h_B : B=h^i(A) \lo X$
    such that
    $f \circ h_B $
    and $g_A\circ (h^i)^{-1}|B$ are ${\cal Y}_{i+4}$-close.
    By (6), $f|B$ and  $g_A\circ (h^i)^{-1}|B$ are
    ${\cal Y}_{i+3}$-close and hence by (2)
    $f \circ h_B$ and $f |B$
    are ${\cal Y}_{i+2}$-close. Then by (4) there is
    a homeomorphism $h_{i+1} : X \lo X$
    such that
    $f\circ h_{i+1}$ and $f$ are ${\cal Y}_{i+1}$-close
    and $h_{i+1} |B=h_B$.
    Note that  the ${\cal Y}_{i+4}$-closeness of $f \circ h_B $
    and $g_A\circ (h^i)^{-1}|B$ is equivalent to
    the ${\cal Y}_{i+4}$-closeness of $f\circ h^{i+1} |A$
    and $g_A$.
     The construction is completed.

    Denote $g_i = f \circ h^i : X \lo Y$. By
    (5) we have that

    (7)  $g_i$ and $g_{i+1}$ are ${\cal Y}_{i+1}$-close\\
    and by (3)

      (8)  $g_j^{-1}({\cal Y}_{i+1})
       \prec_{n-1}g_j^{-1}({\cal Y}_{i})$.
       \\
      Define $g =\lim g_i : X \lo Y$. By (6) we have  $g|A=g_A$ and
       by (2), (7) and (8) we have

     $g^{-1}({\cal Y}_{i+5})
     \prec g_{i+5}^{-1}({\cal Y}_{i+3})
      \prec_{n-1} g_{i+5}^{-1}({\cal Y}_{i+2})
       \prec g^{-1}({\cal Y}_i)$.
      \\
      This implies that $g$ is $UV^{n-1}$ and the proposition
      follows.
         \hfill $\Box$

                    \begin{proposition}
           \label{p2.5}
           Let $A$ be a $\sigma$-$Z$-set in an $n$-dimensional
           \noo space $X$ modeled on a triangulated manifold.
           Then $X\setminus A$ is homeomorphic
            to an $n$-dimensional \noo space
               modeled on a triangulated manifold.
               \end{proposition}
     {\bf  Proof.}
     Represent $A=\cup_{i=1}^\infty A_i$ where $A_i \subset A_{i+1}$
     and $A_i$ is a $Z$-set in $X$.
     Set $X_0=X$, $A_0=\emptyset$ and  $X_i=X \setminus A_i$.
      Note that by Theorem \ref{unknot},
     $X_i$ is homeomorphic to an $n$-dimensional
     \no space modeled on a triangulated manifold. Also note that
     the inclusion of $X_{i+1}$ into $X_i$ is a $UV^{n-1}$-map.
     Fix a complete metric $d$ on $X$.
     Using  Theorem \ref{unknot1} construct for every $i$
     an open cover  ${\cal V}_i$  of $X_i$ and
      a homeomorphism
     $f_i : X_i \lo X_{i+1}$ such that

     (1) mesh${\cal V}_i< 1/2^i$;

    (2)  ${\cal V}_i$  properly
     approaches $A_i$;

     (3)
     $\st{\cal  V}_{i+1} \prec {\cal V}_i$;

     (4) $f_i$ is ${\cal V}_i$-close to the identity map of $X_i$;

     (5)  mesh$(f_0^i)^{-1}({\cal V}_i)< 1/2^i$\\
     where $f_0^0=$id and
     $f_0^{i}= f_{i-1} \circ \dots f_0 : X=X_0 \lo X_{i}$.

      Denote $f =\lim_{i\rightarrow \infty}f_0^i : X \lo X$
      and let us show that $f(X)=X \setminus A$ and $f$ is a
     homeomorphism between $X$ and $X \setminus A$.

     Let $x \in X$. By (3) and  (4), $f(x)$ and $f_0^i(x)$ are
     $ \st^2{\cal V}_i$-close and hence  by (2),
     $f(x)$ is not in $ A_i$. Thus $f(x) \in X \setminus A$.

     Take  $y\in X \setminus A$ and let $x_i \in X_{i}$
     be such that $f_0^i (x_i)=y$. By (3) and (4),
     $f_0^i(x_i)$ and  $f_0^{i+1}(x_{i})$ are ${\cal V}_i$-close
      and hence by (1) and (5),   $d(x_i, x_{i+1})\leq 1/2^i$.
     Denote $x =\lim x_i$. By (3) and (4),
      $f(x_i)$ and  $f_0^i (x_i)$ are $\st^2{\cal V}_{i}$-close
      and hence $f(x)=y$.
       Thus we showed that $f(X)=X \setminus A$.

         Now take   $x_1, x_2 \in X $
         such that $d(x_1, x_2)> 1/2^{i-100}$ for some
         $i > 100$.
         Then by (3) and (5),
           $f_0^i(x_1)$ and  $f_0^i(x_2)$ are not
          $\st^4{\cal V}_i$-close. By (3) and (4),
          $f(x)$ and $f_0^i(x)$ are $\st^2 {\cal V}_i$-close
          for every $x\in X$ and hence $f(x_1)$ and $f(x_2)$ are
          not $\st{\cal V}_i$-close. This implies that $f$ is
          a homeomorphism between $X$ and $X \setminus A$.
              \hfill $\Box$
         \begin{proposition}\label{p2.ad}
         Let $X$ be  an $n$-dimensional \noo space
         modeled on a triangulated manifold.
         Then $X$ can be   split into
         closed subsets
         $X=\cup_{i \in \z} X_i$ such that
         $\{ X_i : i \in \z\}$ is a locally finite cover
         of $X$,
          $X_i\cap X_{j}=\emptyset$ if $|i-j| > 1$,
         $X_i$ and $A_i=X_i \cap X_{i-1}$ are
         homeomorphic to $n$-dimensional \noo
         spaces modeled on triangulated manifolds
         and the sets $A_i$ and $A_{i+1}$ are
         $Z$-sets in $X_i$.
         \end{proposition}
         {\bf Proof.} By Theorems \ref{model} and \ref{unknot}
    we can assume that   $X$ is   the $n$-dimensional
       \no space $R^{2n+2}(n+1)$ modeled on
       the $(2n+2)$-dimensional
        Euclidean space $R^{2n+2}$.
        By a reasoning similar to the one used
        in  the interpretation of \no spaces
        given
         in the introduction (Section 1) of \cite{lev}
        we can represent $X$ as
        $X=R^{2n+2}\setminus K$
        where $K$ is the union of the rational planes
         of $\dim \leq n+1$ (an $m$-dimensional
        plane  of $R^{2n+2}$ is
        rational if it is
        spanned by $m+1$ points with rational coordinates).
    Consider $R^{2n+2}$ as the product
    $R^{2n+2}= R^{2n+1} \times R$
      and   denote by
        $p : R^{2n+2}= R^{2n+1} \times R \lo R $
        and $q : R^{2n+2}= R^{2n+1} \times R \lo R^{2n+1}$
       the projections.

       Choose a discrete sequence of irrational numbers
       $a_i \in R$ indexed by the integers $i \in \z$ such that
       $a_i < a_{i+1}$,
       $\lim_{i\rightarrow -\infty}a_i=-\infty$
       and $\lim_{i\rightarrow \infty}a_i=\infty$.
       Denote
      $X_i=X \cap p^{-1}([a_i,a_{i+1}])$ and
       $A_i=X \cap p^{-1}(a_i)$.
       Since
       $X_i \setminus (A_i \cup A_{i+1})=
       X \cap p^{-1}((a_i,  a_{i+1}))$,
       $X_i \setminus (A_i \cup A_{i+1})$  is
       homemorphic to
        an $n$-dimensional \no space modeled on
       a triangulated manifold.

       Note that  $A_i= p^{-1}(a_i) \setminus K
       \subset  p^{-1}(a_i)\cap
       q^{-1}(R^{2n+1}(n))$ and $ p^{-1}(a_i) \cap K$
       is a countable union of planes of $\dim \leq n$.
       Also note that the intersection of a plane of $\dim \leq n$ in $R^{2n+1}$
       with
        $R^{2n+1}(n)$ is a $Z$-set in  $R^{2n+1}(n)$.
       Thus $A_i$ can be considered as a subspace of
       $R^{2n+1}(n)$ from which countably many $Z$-sets
       are removed.
       Hence
        by Proposition \ref{p2.5},
        $A_i$ is homeomorphic
       to an $n$-dimensional \no space modeled on
       a triangulated manifold. It is easy  to verify
       that $X_i$ is an $n$-dimensional \no space
       and $A_i$ and $A_{i+1}$ are $Z$-sets in $X_i$.
       Then, since  $X_i \setminus (A_i \cup A_{i+1})$
       is homeomorphic
       to an $n$-dimensional \no space modeled on
       a triangulated manifold,  we get
       by Theorem \ref{unknot} that
       $X_i$ is also homeomorphic
       to an $n$-dimensional \no space modeled on
       a triangulated manifold.
       The proposition is proved.
             \hfill $\Box$

         \begin{theorem}\label{resol1}
         Let a complete space $Y$ of $\dim Y \leq n$
         be an absolute extensor in dimension $n$.
         Then there is a $UV^{n-1}$-map $f : X \lo Y$
         from an $n$-dimensional \noo space $X$ modeled
         on a triangulated manifold.
         \end{theorem}
       {\bf  Proof.}
       Split   $X$
      as described in Proposition \ref{p2.ad}.
        For every $i$ choose  a  homeomorphism
         $h_{A_i} : A_i \lo A_0$ such that
         $h_{A_0}=id_{A_0}$.
       Consider a closed embedding  $e : Y \lo A_0$   of $Y$ into
       $A_0$ and  a continuous  retraction $r : A_0 \lo e(Y)$,
       and denote
        $f_{A_i}=e^{-1}\circ r \circ h_{A_i} : A_i \lo Y$.
       Note that $f_{A_i}=f_{A_0} \circ h_{A_i}$.
       By Proposition \ref{p2.3} there is a $UV^{n-1}$-map
       $g_{X_i} : X_i \lo A_i$ such that $g_{X_i} | A_i$ is the identity
       map of $A_i$
       and $g_{X_i}|A_{i+1}=
        h_{A_{i}}^{-1}\circ
        e \circ f_{A_0}\circ h_{A_{i+1}} : A_{i+1} \lo A_i$.
       Set  $f_{X_i}=f_{A_i} \circ g_{X_i} : X_i \lo Y$  and
       note that $f_{X_i}|A_{i+1}=f_{A_{i+1}}$ and
         $f_{X_i}|A_{i}=f_{A_{i}}$. Then the maps $f_{X_i}$
         define the corresponding map $f : X \lo Y$.

         Let us show that $f$ is a $UV^{n-1}$-map.
         Take open
         sets $V $ and $ U$ in $Y$ such that
          $V\subset U$ and the inclusion of
         $V$ into $U$ induces the zero-homomorphism of the homotopy
         groups in $\dim \leq n-1$. Consider a map
         $\psi : S\lo f^{-1}(V)$ from a sphere $S$ of
         $\dim \leq  n-1$ and let us show that
         $\psi$ is null-homotopic inside  $f^{-1}(U)$.
         Since $S$ is compact there are   $i  < j$  such that
         $\psi(S) \subset X_i \cup X_{i+1} \cup \dots \cup X_j$.

         Recall that   $g_{X_j} : X_j \lo A_j$ is a  $UV^{n-1}$-retraction
         onto $A_j \subset X_j$ and
         note that $f^{-1}(V) \cap X_j=
         g_{X_j}^{-1}(W)$ for $ W=f^{-1}_{A_j} (V)$.
         Then by Proposition \ref{p2.1},
         for any map
         $\phi : B \lo X_j$ from a space of $B$ of $\dim \leq n-1$
         such that $\phi(B)\subset f^{-1}(V) \cap X_j$,
           $\phi$ can be homotoped inside $f^{-1}(V) \cap X_j$ and
           relative to $\phi^{-1}(A_j)$ to the map
           $g_{X_j} \circ \phi : B \lo A_j$.

         Thus we can homotope
          $\psi$ inside $f^{-1}(V)$ into a map
         to $X_i \cup \dots \cup X_{j-1}$ and, proceeding by induction
         on $j$, finally homotope $\psi$ inside $f^{-1}(V)$
         into a map to $A_i$.
         Now $\psi$ can be homotoped
         inside $f^{-1}(V)$ to a map to
         $g_{X_{i-1}}\circ \psi : S \lo g_{X_{i-1}}(A_i)=
         h^{-1}_{A_{i-1}}(e(Y))\subset A_{i-1}$.
         Because of the similarity between $A_{i-1}$ and
         $A_0$  we may assume that $A_{i-1}=A_0$. Thus
         without
         loss of generality we may assume that
         $\psi : S \lo e(V)=f^{-1}(V)\cap e(Y) \subset A_0$. Then
       $\psi$ is null-homotopic inside $e(U)\subset f^{-1}(U)$
         and the theorem follows.
           \hfill $\Box$
    \\\\
         Let $f : X \lo Y$ be a map. We say that
         a point $x \in X$ is a regular point of
         $f$ if the collection
         $\{ f^{-1}(U): f(x) \in U, U$ is open in $Y \}$
         is a base of $x$ in $X$.
         Let $A$ be a closed subset of $Y$.  Denote by $X\cup_f A$ the disjoint
         union of $X\setminus f^{-1}(A)$ and $A$, and define
         the topology of $X\cup_f A$ such that
           the topology of   $X\setminus f^{-1}(A)$ is preserved
            and for every $U$ open in $Y$,
         the set $(U \cap A) \cup  f^{-1}(U \setminus A)$
         is open in $X\cup_f A$. Then the space
         $X\cup_f A$ is separable metrizable
         (since so are $X$ and $Y$),  the induced
         function  $f' : X' =X\cup_f A \lo Y$
         defined by $f'(a)=a$ if $a \in A$ and $f'(x)=f(x)$
         if $x \in X\setminus f^{-1}(A)$ is continuous
         and every point of $A$ is a regular point of $f'$.

         A  proof of the following proposition is left to the
         reader.

         \begin{proposition}
         \label{p2.4} Let $f : X \lo Y$ be a $UV^{n-1}$-map
         of $n$-dimensional \noo spaces $X$ and $Y$ and let
         $A \subset Y$ be a $Z$-set in $Y$. Then
         $X'=X\cup_f A$ is an $n$-dimensional \noo space,
         $A$ is a $Z$-set in $X'$
         and the induced map $f' : X' \lo Y$ is $UV^{n-1}$.
         \end{proposition}

         \begin{theorem}
         \label{resol2}
         Let $f : X \lo Y$ be a $UV^{n-1}$-map from
          an $n$-dimensional
           \noo space $X$ modeled on a triangulated manifold
           to an $n$-dimensional \noo space $Y$. Then there is
          a $UV^{n-1}$-map $g : X \lo Y$  for which there is
          a $\sigma$-$Z$-set $A \subset Y$ such that
          $g^{-1}(A)$ is a $\sigma$-$Z$-set in $X$ and
            $g^{-1}(Y\setminus A)$ consists of
            regular points of $g$.

           \end{theorem}
    {\bf  Proof.}
    Fix  complete metrics in $X$ and  $Y$.
    Take a sequence  of open covers ${\cal X}_i, i=1,2\dots$
    and  $Z$-embeddings
    $\psi_i : X \lo X, i=2,3,\dots$ such that mesh${\cal X}_i \leq 1/2^i$,
    $\psi_i$ is ${\cal X}_i$-close to the identity map
    and the $Z$-sets $A_i =\psi_i(X)$ are pair-wise  disjoint.

   Set $f_1=f$, $A_1=\emptyset$ and let ${\cal Y}_1$ be an open cover of $Y$
   with  mesh${\cal Y}_1\leq 1/2$.
     We are going to
    construct
    for every $i$ an open cover ${\cal Y}_i$  and
    a $UV^{n-1}$-map $f_i : X \lo Y$ such that:

    (1)
    mesh${\cal Y}_i\leq 1/2^i$ and
   $\st^5{\cal Y}_{i+1}\prec_{n-1} {\cal Y}_i$;

    (2)
     $f_i$ and $f_{i+1}$ are $ {\cal Y}_{i}$-close;

     (3) $f_i^{-1}(f_i(A_i))=A_i$, $f_i|A_i$ is a $Z$-embedding
     of $A_i$ in $Y$ and
     $f_i|A_j= f_j|A_j$ for $i> j$;

     (4) every point of $A_1\cup \dots\cup A_i$ is
     a regular point of $f_i$  and
     the family
     ${\cal A}_i=\{ f_i^{-1}(U) : f_i^{-1}(U) \cap (A_1 \cup\dots\cup A_i)
     \neq  \emptyset,  U\in {\cal Y}_i \}$ refines ${\cal X}_i$.\\

     Assume that the construction for the  indices$\leq i$ is
     completed. Proceed to $i+1$ as follows.

     Let ${\cal U}={\cal Y}_i$.
      By Theorem \ref{unknot1}
     there is an open cover ${\cal V}$ of $Y$ such that
     the conclusions of Theorem \ref{unknot1} hold for $f$
     replaced by
     the map  $f_i: X \lo Y$.

     Approximate $f_i$ by a $Z$-embedding $\phi : X \lo  Y$ such
     that  $\phi(X)$ does not intersect
     $f_i(A_1 \cup \dots \cup     A_i)$ and $\phi$ is
     ${\cal V}$-close to $f_i$.
     Denote $B=\phi(A_{i+1})$ and $X'=X \cup_{f_i} B$, and
     let $f' : X' \lo Y$ be the map induced by $f_i$.
     By Proposition \ref{p2.4} and Theorem \ref{unknot}
     the space $X'$ is homeomorphic to an $n$-dimensional
     \no space modeled on a triangulated manifold.
     Then by Theorem \ref{unknot1} one can choose
     a homeomorphism $h' : X \lo X'$ such that
     $f' \circ h'$ is arbitrarily close to $f_i$ and
     $h'(a)=a $ for $a \in A_1 \cup \dots\cup A_i$. In
     particular $h'$ can be chosen so that
      $f' \circ h'$ is  ${\cal V}$-close
     to $f_i$.
     Then again by Theorem \ref{unknot1} there is
    a  homeomorphism  $h : X \lo X$ such that
     $h(a)=(h')^{-1}(\phi(a))$ for $a \in A_{i+1}$,
     $h(a)=a$ for $a \in A_1 \cup \dots \cup A_i$
     and
     $f_{i+1}= f' \circ {h'}\circ h $ is
     $\cal U$-close to $f_i$.
     Note that every point of $A_1 \cup \dots\cup A_{i+1}$
     is a regular point of $f_{i+1}$.

     Choose
      an open cover ${\cal Y}_{i+1}$
     of $Y$ such that the properties (1),
     (2) and (4) are
     satisfied. The construction is completed.

     Define $g=\lim f_i : X \lo Y$. The properties
     (1) and (2)
      imply that

     (5)  $g^{-1}({\cal Y}_{i+3})
     \prec f_{i+3}^{-1}({\cal Y}_{i+2})
     \prec f_{i+3}^{-1}( {\cal Y}_{i+1})
      \prec g^{-1}({\cal Y}_i) $.
       \\
      Since ${\cal Y}_{i+2} \prec_{n-1}{\cal Y}_{i+1}$
      and $f_{i+3}$ is $UV^{n-1}$ we have
      $f_{i+3}^{-1}({\cal Y}_{i+2})
     \prec_{n-1} f_{i+3}^{-1}( {\cal Y}_{i+1})$ and
        therefore $g$ is  $UV^{n-1}$.
     The properties (3-5)  imply that for every $i$,
     $g|(A_1 \cup \dots \cup A_i)$ is a $Z$-embedding of
      $A_1 \cup \dots \cup A_i$ into $Y$ and
      every point of $A_1 \cup \dots \cup A_i$ is a  regular
      point
      of $g$.  Denote by $C$ the set of all regular points of $g$.
      Then
       $C$ is $G_\delta$ in $X$ and since  $A_i \subset C$ for every
      $i$ we have that $X \setminus C$ is a $\sigma$-$Z$-set in
      $X$. Hence by Proposition \ref{p2.1} one can choose
       a map $\psi : Y \lo X$ such that $\psi(Y) \subset X
       \setminus C$ and
         $ g\circ \psi $ is arbitrarily close to
       the identity map  $Y$. Thus
             $Y \setminus g(C)$ is
           a $\sigma$-$Z$-set in
      $Y$ and
           the conclusions of
      the theorem  hold for $A =Y \setminus g(C)$.
       \hfill $\Box$
        \\\\
        {\bf  Proof of Theorem \ref{resol}.}
        By Theorems \ref{resol1} and \ref{resol2}
         there is a $UV^{n-1}$-map
        $g : X' \lo Y$ from an $n$-dimensional \no
        space $X'$ modeled on a triangulated manifold
        such that there is a $\sigma$-$Z$-set $A$ in $Y$
        for which $g^{-1}(A)$ is a $\sigma$-$Z$-set in $X'$
        and $g$ embeds $X=X' \setminus g^{-1}(A)$ into $Y$.
        By Proposition \ref{p2.5}, $X$ is homeomorphic
        to an $n$-dimensional \no space and the theorem follows.
         \hfill $\Box$

    \end{section}

    \begin{section}{Auxiliary constructions and properties}\label{aux}

      \subsection{A few general properties}\label{general}
    Let $M$ be  an
    $m$-dimensional  triangulated manifold $M$ with
    $m\geq 2n+1$. Set $k=m-n-1$. Fix a triangulation  $\cal T$
    of $M$ and
     embed  $M$ into
    a Hilbert  space $H$ by a map which is linear
    on each simplex of ${\cal T}$.
      Let $K$ be a countable union of  planes in
      $H$
     of $\dim \leq k$ and denote $X=M\setminus K$
     (it is shown
     in Introduction of \cite{lev}
      that $M(k)$ admits such a representation). It is clear
       that $ M \cap K$ can  be represented as
      a countable union of simplexes of $\dim \leq k$
      PL-embedded in $M$.

      Let $f : B^q \lo \Int M$ be a PL-embedding of
      a $q$-dimensional ball $B^q$ with $q \leq n$.
      Then $f$ can be arbitrarily closely
      approximated by a PL-embedding $f' : B^q \lo \Int M$
     such that $f'(B^q)\subset X$. To show that
      identify $B^q$ with $f(B^q)$ and extend
      the embedding $B^q \subset \Int M$ to
      a PL-embedding $B^m=B^q \times B^{m-q}$
      of an $m$-dimensional ball
      $B^m$ with $B^q=B^q \times O$.
      Since $K$ is a countable union of simplexes
      of $\dim \leq k$ PL-embedded in $M$ one can choose
      a point $a \in B^{m-q}$ arbitrarily close to $O$
      such that $B^q \times a \subset X$. This way we get
      the required approximation of $f$.

   Let  ${\cal P}$ be a decomposition of
    of $M$. Then
     there is an open subset $ M'$ of $M$ containing
    $X$ such that the finite intersections of
    ${\cal P}$ restricted to $ M'$ have dense subsets
    lying in $X$. Indeed, take  a subdivision  ${\cal T}'$
    of $\cal T$ such that ${\cal T}'$ underlies $\cal P$.
    Then
    one can  easily verify that for every simplex $\Delta$
    of ${\cal T}'$ having an open subset lying
    outside $X$, the entire simplex $\Delta$ lies
    outside $X$. Remove from $M$ all the simplexes
    of ${\cal T}'$ lying outside $X$ and get the open subset
    $ M'$ with the required properties.

    Let  $F$ be a PL-subcomplex of $M$ such that
   $\cal P$  forms
    a partition on $M\setminus F$.
     Then for every open subset
    $ M'$ of $M$ containing $X$  we have that
    for every finite intersection $P$
    of ${\cal P}$ the inclusion
    $(P\cap M') \setminus F \subset P \setminus F$
    induces an isomorphism of the homotopy
    groups in co-dimensions $\geq m-n+1$
    (=in dimensions$ \leq \dim(P\setminus F) -(m-n+1)$).
      Indeed,
     let $P$ be a finite intersection of
     ${\cal P}$ with $t=\dim P \setminus F \geq m-n+1$.
     Since $M \setminus X$ as a countable union
     of simplexes of $\dim \leq k$ PL-embedded in $M$,
    we have that $(P\setminus F)\setminus X$
    is also a countable union of simplexes of $\dim \leq k$
     PL-embedded in $P \setminus F$.
      Then
        every map of a sphere of $\dim \leq t-(m-n+1)$
     into $P \setminus F$ can be homotoped into
     $(P \setminus F) \cap X \subset
     (P \setminus F) \cap M'$.
     Now  let
     $f : S \lo ({ M'}\cap P) \setminus F$ be
     a map from a sphere of $\dim \leq t-(m-n+1)$.
       Since
       $(P\setminus F)\setminus X$
    is  a countable union of simplexes of $\dim \leq k$
     PL-embedded in $P \setminus F$ we have that
      $P\setminus(X \cup F \cup f(S))$
     is also a countable union of simplexes of $\dim \leq k$
     PL-embedded in $P \setminus F$.
     Then if
      $f$ is null-homotopic
     in $P \setminus F$ we get that $f$ is null-homotopic
     in $(P \setminus F)\cap (X \cup f(S))$.
     Thus $f$ is null-homotopic in
     $({ M'}\cap P) \setminus F$ and we are done.

    In particular we get that  if
  $\cal P$  forms
    an $l$-co-connected partition on $M\setminus F$
    with $l \geq m-n+1$
    then
     $\cal P$ restricted to ${ M'} \setminus F$
     is also $l$-co-connected. Note that the correspondence
     $\mu$ sending $P \in {\cal P}$ to
     $P \cap M'$ induces an $n$-matching between
     ${\cal P}$ restricted to $M \setminus F$
     and ${\cal P}$ restricted to $M'\setminus F$, see
     \ref{l-match}.

     Let $X=M(k)$. We leave to the reader to verify
     that  for every PL-subcomplex $L$ of $M$ with
     $\dim L \leq n$, $X \cap L$ is a $Z$-set in $X$.

     \subsection{  Creating intersections}\label{create}
    Let $M$ be a triangulated $(l-1)$-co-connected
    manifold such that $m=\dim M \geq 2q+1, q=m-l+2$.
    Assume that
    $F$ is a PL-subcomplex with $\dim F \leq l-2$,  $\cal P$
    is a decomposition
    of $M$ such that $\cal P$ forms an $(l-1)$-co-connected
    partition on $M \setminus F$.
    Consider $t+1$ distinct
    elements  $P_0,\dots, P_t$ of $\cal P$, $t=m-l+2$,
    such that
    their intersection is empty but
    the intersection of any $t$  of them is not
    empty.

     Let us show how  the construction  2.5 of \cite{lev}
    can be used to create an intersection of $P_0, \dots, P_t$.
    Use the notation of 2.5 of \cite{lev} and assume that
    $S_P= \emptyset$. Then for every
    simplex $\Delta'$ of $\Delta$, $S_P * \Delta'=\Delta'$
    and we start the construction for $0$-dimensional
    simplexes
    $\Delta'$ with any $f_{\Delta'} : \Delta'  \lo
   \Int ( P(\Delta')\cap U)$.
    Now we can follow  2.5 of \cite{lev} to
   construct $f_{\partial \Delta}$ and after that
    to modify
    ${\cal P}$ to ${\cal P}^\pi$.  It is easy to see   that
    the modifications
    of $P_0, \dots, P_t$ will intersect on $M \setminus F^\pi$.

    Now suppose that we have a (possibly countable)
    collection
    ${\cal C}=\{ (P_0, \dots , P_t) : P_i \in {\cal P}\}$
    of $(t+1)$-tuples of elements of ${\cal P}$
    of the type described above and we need to create
    an intersection for every $(t+1)$-tuple from $\cal C$.
     Assume that ${\cal W}$
    is an open cover of $M$ such that
    for every
    $(t+1)$-tuple $(P_0, \dots, P_t) \in {\cal C}$
     there is a set $W \in \cal W$
    such that $P_0 \cup \dots \cup P_t \subset W$
    and the inclusion $P_0 \cup \dots \cup P_t \subset W$
    induces the zero homomorphism
    of the homotopy groups in $\dim \leq m-l+1$.

    Then  the construction  2.7 of \cite{lev}
    applies to modify $M$ to an open subset $M'$ of
    $M$, $F$ to a PL-subcomplex of $F'$ of $M'$,
    each element $P$ of $\cal P $ to a PL-subcomplex
    $P'$ of $M'$ such that the following properties are satisfied:
    $\dim M\setminus M' \leq q$, $\dim F' \leq l-2$,
    ${\cal P}' =\{ P' : P \in {\cal P} \}$ is
    a decomposition of $M'$ which forms
    an $(l-1)$-co-connected partition on $M' \setminus F'$,
    $P' \subset \st (P , \st{\cal W})$ for each $P \in \cal P$,
    all the finite intersections of ${\cal P}$ restricted to $M \setminus F$
    are preserved in ${\cal P}'$ restricted to $M' \setminus F'$
     (via the natural correspondence between ${\cal P}$ and ${\cal
     P}'$)
    and the only new finite intersections of ${\cal P}'$
     that are created on
    $M'\setminus F'$ are for the $(t+1)$-tuples in $\cal C$.

       \subsection{  A $t$-matching of partitions}\label{l-match}
       Let $M_1$ and $M_2$ be
       triangulated manifolds and
       let ${\cal P}_1$ and ${\cal P}_2$ be partitions
           of $M_1$ and $M_2$ respectively.
       A  one-to-one correspondence
       $\mu : {\cal P}_1 \lo {\cal P}_2$
       is said to be
       a $t$-matching
       if for any  $P_0, P_1,\dots P_t\in {\cal P}_1$  we have that
       $P_0 \cap \dots \cap P_t \neq \emptyset$ if and only
       if $ \mu(P_0) \cap \dots \cap \mu(P_t) \neq  \emptyset$.
       In other words, $\mu$ is a $t$-matching if it preserves
       the intersections in co-dimensions$\leq t$.
       Thus $\mu$ is a matching if it is a $t$-matching for
        every $t=0,1,2, \dots$.

        Denote by $F_i$ the union of all finite intersections of
        ${\cal P}_i$ of $\dim \leq \dim M_i -t-1$.
        Then $\mu $ induces a matching of partitions when
        ${\cal P}_1$ and ${\cal P}_2$ are restricted
        to $M_1 \setminus F_1$ and $M_2 \setminus F_2$
        respectively if and only if $\mu$ is a $t$-matching
        between ${\cal P}_1$ and ${\cal P}_2$.

    \subsection{  Improving the level of matching
    of  partitions}\label{level}
     Let $M_i,i=1,2$ be $(l_i-1)$-co-connected
       triangulated manifolds such that
          $m_i=\dim M_i  \geq 2(m_i -l_i +2) +1$
          and $m_1 -l_1=m_2-l_2$.
          Denote $t=m_i-l_i+1$.
          Assume that
           ${\cal P}_1$ is an $(l_1-1)$-co-connected
        partition of $M_1$,
        $F_2$ is a  PL-subcomplex
        of $M_2$, $\dim F_2 \leq l_2-2$,
        ${\cal P}_2$ is a decomposition of $M_2$ forming
         an $l_2$-co-connected
         partition
        of $M_2 \setminus F_2$ and
        $\mu : {\cal P}_1 \lo {\cal P}_2$ is
           a one-to-one correspondence such that $\mu$ induces
           a $t$-matching between ${\cal P}_1$ and ${\cal P}_2$
            restricted to $M_2 \setminus F_2$.
            Let us show how
           to turn $\mu$ into a $(t+1)$-matching simultaneously
           with improving connectivity of ${\cal P}_2$.

           By $M'_2$, $F'_2$ and ${\cal P}'_2$ we  denote
          the modifications of $M_2$, $F_2$ and ${\cal P}_2$
           and we always assume that $M'_2$ is an open subset of
           $M_2$, $M_2\setminus M'_2$ is PL-presented in $M_2$,
           $F'_2$ is a PL-subcomplex of $M'_2$, ${\cal P}'_2$
           is a decomposition of $M'_2$ which forms a partition
           on $F'_2 \setminus M'_2$.
            By $\mu': {\cal P}_1 \lo {\cal P}'_2$ we denote
            the correspondence sending $P \in {\cal P}_1$
            to the modification of $\mu(P)$ in ${\cal P}'_2$.

           Apply the construction 2.8 (improving
           the total connectivity of a partition) of \cite{lev}
           to modify $M_2$, $F_2$ and ${\cal P}_2$ to
           $M'_2$, $F'_2$ and ${\cal P}'_2$ such that
           $\dim M_2\setminus M'_2 \leq m-l_2$,
           $\dim F'_2\leq l_2 -2$ and ${\cal P}'_2$  is
           $(l_2-1)$-co-connected on ${ M}'_2 \setminus F'_2$.\\
\\
           The next step is to
              create
             the missing  intersections of ${\cal P}'_2$ in
            $\dim =m-t-1$.
           Denote

           ${\cal C}=\{ (P_0,\dots P_{t+1}): P_0, \dots , P_{t+1}$
           are distinct elements in ${\cal P}'_2$ such that

           $P_0 \cap \dots\cap P_{t+1} \cap (M'_2 \setminus F'_2)=\emptyset$
           and ${\mu'}^{-1}(P_0) \cap \dots \cap {\mu'}^{-1}(P_{t+1}) \neq
           \emptyset\}$.\\
           Apply  \ref{create} to modify
           $M'_2$ , $F'_2$ and ${\cal P}'_2$
             to create intersections for the $(t+2)$-tuples
            in $\cal C$.\\

           Let us make   $\mu'$ induce a matching.
           We will  say that a  non-empty finite
           intersection $P=P_0\cap \dots\cap P_s$
  of distinct elements in ${\cal P}_1$ is brought from ${\cal
  P}'_2$  if
  $({\mu'}(P_0) \cap \dots \cap{ \mu'}(P_s))
  \cap(M'_2 \setminus F'_2) \neq \emptyset$.
  Similarly we say that
  a   finite intersection $P=P_0\cap \dots\cap P_s$
  of distinct elements in ${\cal P}_2$ is brought from a finite
  intersection of  ${\cal
  P}_1$   if $(P_0 \cap \dots \cap P_s )\cap(M'_2 \setminus F'_2)
  \neq \emptyset$
  and
  ${\mu'}^{-1}(P_0) \cap \dots \cap {\mu'}^{-1}(P_s) \neq \emptyset$.

           Denote

           $F_1 = $the union of the finite intersections of ${\cal P}_1$ of
           $\dim < m-t-1$ which are not brought from ${\cal P}_2$;

            $F^+_2 =$the union of the finite intersections of
            ${\cal P}_2$ which are not brought from
            the finite intersections of ${\cal P}_1$ of $\dim <
            m-t-1$

           ${F}^{++}_2=$the union of finite intersections of
           ${\cal P}_2$ which are not brought from
            the finite intersections of ${\cal P}_1$ of $\dim =
            m-t-1$

       Since
            $F^+_2$ and $F^{++}_2$ are unions of finite
            intersections of ${\cal P}_2$  then
            for every finite intersection $P$ of ${\cal P}'_2$
            such that
            $(P   \setminus F'_2) \setminus ( F^+_2 \cup
            F^{++}_2) \neq \emptyset$ we have
            $(P \setminus F'_2) \cap (F^+_2 \cup F^{++}_2)
            \subset \partial (P \setminus F'_2)$,
            see 2.2 of \cite{lev}.
            Hence replacing  $F'_2$ by the union
            $F'_2 \cup F^+_2 \cup F^{++}_2$ we get that
            ${\cal    P}'_2$ remains to be an $(l_2-1)$-co-connected
            partition on $M'_2\setminus F'_2$, $\dim F'_2 \leq
            l_2-2$ and every finite intersection of ${\cal P}'_2$
            (on $M'_2\setminus F'_2$) is brought from ${\cal P}_1$
            by the correspondence $\mu'$.
            Similarly we conclude that ${\cal P}_1$ restricted to
            $M_1\setminus F_1$ is
            $(l_1-1)$-connected.

            Thus $\mu'$ becomes
            a matching of partitions when
            ${\cal P}_1$ and ${\cal P}'_2$ are
            restricted to
            $M_1 \setminus F_1$ and $M'_2\setminus F'_2$
            respectively.
            Since
             $\dim F_1 \leq m_1 -t-2= l_1-3$ and $M_1$ is
             $(l_1-1)$-co-connected  we get that $M_1\setminus
             F_1$ is also  $(l_1-1)$-co-connected. Then by
             2.3 of \cite{lev} we get that $M'_2 \setminus F'_2$
             is $(l_2-1)$-co-connected.
             Now we can  apply the construction 2.9 (absorbing simplexes)
             of \cite{lev} to modify
             $M'_2$, $F'_2$ and ${\cal P}'_2$ in order to reduce
             the dimension of $F'_2$ to $\dim \leq l_2 -3$
             leaving the other characteristics of
             $M'_2$, $F'_2$ and ${\cal P}'_2$ unchanged.

            Thus we finally get that $\mu'$ is a $(t+1)$-matching between
            ${\cal P}_1$ and ${\cal P}'_2$ restricted to
            $M'_2\setminus F'_2$,
            $\dim M_2\setminus M'_2 \leq m_2 -l_2$,
             $\dim F'_2 \leq l_2-3$ and ${\cal P}'_2$ is
             $(l_2 -1)$-connected on $M'_2\setminus F'_2$.
             \\

             The procedure described above can be  used
             iteratively as follows.
             Assume  that $\dim M_i \geq 2n+1$,
             $M_i$ is $(n-1)$-connected(=$(m_i-n+1)$-co-connected), ${\cal P}_1$ is
             an $(m_1-n+1)$-co-connected
             partition of $M_1$ and
             ${\cal P}_2$ is a partition of $M_2$ which
             admits a $0$-matching
             $\mu: {\cal P}_1 \lo {\cal P}_2$. Then
             repeating inductively the above procedure
             we can modify
             $M_2$, $F_2$,  ${\cal P}_2$ and $\mu$
             to $M'_2$, $F'_2$, ${\cal P}'_2$
             and $\mu' : {\cal P}_1 \lo {\cal P}'_2$
             such that $\dim M_2 \setminus M'_2 \leq n$,
             $\dim F'_2 \leq m_2 -n-1$, ${\cal P}'_2$ is
             $(m_2-n+1)$-co-connected on $M'_2 \setminus F'_2$
              and $\mu'$
             induces an $n$-matching between ${\cal P}_1$
             and ${\cal P}'_2$ restricted to
             $M'_2 \setminus  F'_2$.

   \subsection{ A remark on improving connectivity of intersections}
   \label{improv}
   In the construction 2.4 of  \cite{lev} we consider
    a decomposition $\cal P$ of a triangulated
    $m$-dimensional manifold
    $M$ and a PL-subcomplex $F$ such that
    ${\cal P}$ forms  a partition on $M \setminus F$. We fix
    a finite intersection $P$ of $\cal P$ whose connectivity
    has to be improved, take
      a PL-embedding  $f $ of
   a sphere  $S_P$  into $\Int(P \setminus F)$,  extend $f$
   to a PL-embedding  $f_{\partial \Delta} $
   of a larger $(q-1)$-dimensional sphere
   $S^{q-1}$, identify $f_{\partial \Delta} (S^{q-1})$
   with the boundary $\partial B^q$ of a $q$-dimensional ball
   $B^q$
   and finally (using 4.2 of \cite{lev}) observe
   that the  PL-embedding of $\partial B^q$ can be extended
    to a PL-embedding  of an $m$-dimensional ball
    $ B^m =B^q \times B^{m-q} \subset \Int M$ such that
    $B^m \cap F =(B^q \cap F)\times B^{m-q}$ and
    for every  $P'\in\cal P$,
    $B^m \cap P' =(B^q \cap P')\times B^{m-q}$.

    Now  it is clear
    that      the original embedding of $\partial B^q=\partial B^q
    \times O$ can be replaced by any embedding
    $\partial B^q \times a$, $a \in \Int B^{m-q}$.
    Let $X =M(k)$ and $m-q >k$. Then,
    by \ref{general}, $M\setminus X$ is a countable union of
    simplexes of $\dim \leq k$ PL-embedded in $M$
    and hence $a \in B^{m-q}$ can be chosen  to be arbitrarily
    close to $O$ and such that $ B^q  \times a \subset X$.
    Thus without loss of generality we can replace
    the original embedding $f_{\partial \Delta}$
    by an arbitrarily close  embedding whose
    image is contained in $X$.

       \subsection{ A remark on absorbing simplexes}\label{rem-absorb}
       Let $M_i$, $i=1,2$ be an $l_i$-co-connected $m_i$-dimensional
        manifold such that $m_1-l_1=m_2-l_2\geq 0$  and
        $m_i \geq 2(m_i-l_i +1)+1$. Let $F_i$ be a PL-subcomplex of
        $M_i$ with $\dim F_i \leq l_i-1$ and let  ${\cal P}_i$ be
        a decomposition of $M_i$ such that ${\cal P}_i$ forms
        an $l_i$-co-connected partition on $M_i \setminus F_i$
        and there is a one-to-one correspondence
        $\mu : {\cal P}_1 \lo  {\cal P}_2$ which induces
        a matching of partitions
        when ${\cal P}_1$ and ${\cal P}_2$ are
        restricted  to $M_1 \setminus F_1$ and
        $M_2\setminus F_2$ respectively.

              We want to reduce the dimension of $F_2$ to $l_2-2$.
         Following the construction 2.9 (absorbing simplexes)
        of \cite{lev}
        fix a sufficiently fine triangulation of $M_2$
        underlying ${\cal P}_2$ and $F_2$.
        Assume such that ${\cal A}$ is a cover $M_1$ such that the
        elements of $\cal A$ are unions of elements of ${\cal
        P}_1$ and for every simplex $\Delta \subset F_2$
        of $\dim=l_2-1$ there is an element $A \in {\cal A}$
        such that $\mu^{-1}(\st(\Delta, {\cal P}_2)) \subset
        A$ and the inclusion
        $\mu^{-1}(\st(\Delta, {\cal P}_2))\setminus F_1
        \subset   A\setminus F_1$
        induces the zero-homomorphism
        of the homotopy groups in dimensions$\leq m_1-l_1$
        (=co-dimensions$\geq l_1$). Then
        the inclusion
        $\st(\Delta, {\cal P}_2)\setminus F_2
        \subset  \mu( A)\setminus F_2$  induces
        the zero-homomorphism of the homotopy
        groups  in $\dim \leq m_2-l_2$, see 2.3 of \cite{lev}.
        Now it is easy to derive from 2.9 of \cite{lev}
        that the construction of absorbing the
        $(l_2-1)$-dimensional simplexes of $F_2$ can be carried
        out such that the modification ${\cal P}'_2$ of ${\cal P}_2$
         will
        refine $\mu(\st^2{\cal A})$ and for every $P\in {\cal P}_2$,
        the modification of $P$ will be contained
        in $\st (P, \mu(\st^2{\cal A}))$.

        It is also easy to derive from 2.9 of \cite{lev} that
        if $K$ is a  PL-subcomplex of $M_2$ such that
        $\dim K \leq l_2-2$ and $K \cap F_2 =\emptyset$
        then the modification ${\cal P'}_2$ of ${\cal P}_2$
        can be constructed so that  $K$ will be contained
        in the modification $M'_2$ of $M_2$ and ${\cal P}_2$
        will be left unmodified on a small neighborhood
        of $K$.

    \subsection{ A remark on digging holes for improving connectivity}
    \label{rem-dig}
   In the construction 4.2 (digging holes for improving
    connectivity) of \cite{lev} we extend a PL-embedding
    of an $(q-1)$-dimensional sphere $S^{q-1}$ into
    a manifold $M$ to a PL-embedding of  a ball
    $B^m=B^q \times B^{m-q}$ such that
    the ball $B^q$ bounds $S^{q-1}$ and $B^m$ has
    special properties with respect to a decomposition
    ${\cal P}$ of $M$.
    It is  noted  in the end of 4.2 of \cite{lev} that if $S^{q-1}$
    can be contracted to a point in an open subset $W$ of $M$
    then $B^m$ can be chosen to  be in $W$. We will need a slightly
    improved  control on $B^m$. Namely,
    if $B^q_{\#}$ is any $q$-dimensional ball with
    $\partial B^q_{\#}=S^{q-1}$ and  $\phi_\#: B^q_{\#}  \lo M$ is a map
    which is the identity map on $S^{q-1}$ then the ball $B^q$ can be
    chosen to be arbitrarily close to $\phi_\#(B^q_{\#})$ in the sense
    that $B^q$ can be chosen to be  $B^q=\phi(B^q_\#)$ for
    a sufficiently  close approximation of  $\phi_\#$ by
    a PL-embedding $\phi :  B^q_{\#}  \lo M$ such that
    $\phi$ is the identity map on $S^{q-1}$. This can be done as
    follows.

    The case
      $m=3$ and $q=1$ can be visualized directly.
     Assume that $m\geq 4$ (actually we can always assume
     that $m \geq 4$  restricting ourselves to
      \no spaces modeled
     on manifolds of $\dim \geq 4$).
     The space $R^m$ in the beginning of 4.2 of \cite{lev}
     can be chosen so that $\phi_\#(B^q_\#) \subset R^m$.
     Then after replacing the original embedding of
      $R^m$ by $e : R^m \lo e(R^m)= R^m \subset M$
     we also have $\phi_\#(B^q_\#) \subset R^m$.

     Approximate $\phi_\#$ by a PL-embedding
     $\phi'_\#: B^q_\# \lo R^m$ such that
     $\phi'_\#$
     coincides with $\phi_\#$ on
     $\partial B^q_\#$ and $\phi'_\#$
     sends a small
     neighborhood of $\partial B^q_\#$ into $B^q$.
     Then  the block bundle $\eta$
     can be chosen to be
     so close to $S^{q-1}$ that
      $E(\eta)\cap \phi'_\#(B^q_\#)=E(\eta)\cap B^q$.
         Since $m \geq 4$, it follows from
    Unknotting Theorem (see
     Theorem 10.1 of \cite{h}) that
    there is
    a PL-homeomorphism $e': R^m \lo R^m$ such that
    $e'$ does not move the points of $E(\eta)$ and
    $e'$ sends $\phi'_\#(B^q_\#)$ onto $B^q$.

    Now replace the embedding $e$ by
    $e'\circ e : R^m \lo R^m \subset M$.
    Following the construction
   4.2 of \cite{lev} we can  choose the point $a \in B_1^{m-q}$
   to be arbitrarily close to the center $O$ of  $B_1^{m-q}$.
   Then   the homeomorphism $\psi : M \lo M$ can be chosen
   to be  arbitrarily close
   to the identity map. Thus  the final embedding of $B^q$ can
   be constructed to be arbitrarily close
      to $\phi'_\#(B_\#^q)$ and hence it can be made
      close to $\phi_\#(B_\#^q)$.

   \subsection{  A remark on moving to a rational position}
   \label{rem-mov}
   The reasoning of construction 2.11 (moving to a rational position) of
   \cite{lev} applies to show the following property.
   Assume $\cal P$ is a decomposition
   of a triangulated manifold $M$ which is represented
   as  the union $M =U \cup V$ of two open subsets $U$ and $V$
   such $\cal P$ restricted to $V$ is a rational decomposition of
   $V$. Then there is  a PL-homeomorphism $h : M \lo M$
   such that $h({\cal P})$ is a rational decomposition,
   $h(x)=x$ for every $x \in X \setminus U$ and $h$ can be chosen
   to be arbitrarily close to the identity map.

   Indeed, let  a rational triangulation $\cal T$ of $M$ be
   such that
   $\st(U \setminus V, {\cal T})\cap
    \st (V\setminus U,{\cal T})=\emptyset$.
    Embed $M$ into a Hilbert space $H$ by a map which is linear
    on the simplexes of ${\cal T}$.
    Take a subdivision
    ${\cal T}'$ of ${\cal T}$ such that
     the simplexes of ${\cal T}'$ are linear in $H$,
    ${\cal T}'$ underlies
    ${\cal P}$ and  the simplexes of ${\cal T}'$ contained in
    $ \st (V\setminus U,{\cal T})$ are rational.
    Similarly to 2.11 of \cite{lev}, we can define a map $h'$
     from the vertices of ${\cal T}'$ to rational
     points in $M$ such that $h'$ does not move the vertices lying
      $ \st (V\setminus U,{\cal T})$ and $h'$ slightly moves
      the other vertices such that
       each of them will not leave every  simplex
     of ${\cal T}$
     to which it belongs. Then $h'$ can be linearly
     extended  a PL-map $h: M \lo M$ and
     $h$ will be a homeomorphism if $h'(v)$ and $v$ are sufficiently close
     for every vertex $v$ of ${\cal T}'$.

      \subsection{  A remark on the discretization
      of maps' images  }\label{rem-discr}
      Assume that in the construction 2.6 of
      \cite{lev} the finite intersections
      of
       the decomposition $\cal P$ have dense subsets on
       $X=M(k)$ where $k=\dim M -r -1$. Then the PL-subcomplex
        $R$ of $M$ can be chosen to be rational.
        Indeed, using \ref{rem-mov}
        the set $R$
       can be moved to a rational subcomplex with
       the needed properties. Note that,  since
       the finite intersections of $\cal P$ are dense
       on $X$,
       every rationally presented subset of $M$ of
       $\dim \leq k$ is
       nowhere dense on the finite intersections of ${\cal P}$.

       This implies that  if in the constructions 2.7 and 2.9
       of \cite{lev}
       the decomposition ${\cal P}$ is rational such that
       the finite intersections of ${\cal P}$ have
       dense subsets on $X=M(n)$ with $\dim M=m\geq 2n+1$ and
       $n \geq  q$
        then in 2.9 of \cite{lev}
       the sets $R^{\partial \Delta}$ and $R$ can be chosen
       to be rationally presented  and in 2.9 of \cite{lev}
       the set $R$ can be chosen to be a rational subcomplex.

    \subsection{  A subdivision of a partition  }\label{subdiv}
    Let ${\cal P}$ be a partition of  a triangulated
    manifold $M$. A partition ${\cal P}_\#$ of $M$ is said to be
    a subdivision of ${\cal P}$ if ${\cal P}_\#$ refines
    ${\cal P}$.

    For every open cover $\cal V$ of $M$ there is a subdivision
    ${\cal P}_\#$ of ${\cal P}$  such that ${\cal P}_\#$ refines
    $\cal V$ and every finite intersections of ${\cal P}_\#$
    is contractible. Indeed, arrange the elements of ${\cal P}$
    into a sequence $P^1, P^2, \dots$. Take a triangulation
    $\cal T$ of $M$ such that $\cal T$ underlies $\cal P$ and
    $\st \cal T$ refines $\cal V$. Denote by  ${\cal T}_i$
    the $i$-th barycentric subdivision of ${\cal T}$.
    Define ${\cal P}^i_\#$ as the collection
    $\{ \st(v, {\cal T}_{i+1})\cap P^i: v$ is a vertex of
    ${\cal T}_i$ and $v \in P^i \}$. Then
    ${\cal P}_\# =\cup_i {\cal P}^i_\#$ is a partition of $M$
    having the required properties. Note that ${\cal P}_\#$
    will be  a rational partition if ${\cal P}$ is
    a rational partition and we choose ${\cal T}$ to
    be a rational triangulation.

    Assume that $\cal P$ and ${\cal P}_\#$ are $l$-co-connected
    partitions of $M$ such that ${\cal P}_\#$ is
   a subdivision of $\cal P$.
    Assume that each $P_\# \in {\cal P}_\#$ is modified to $P'_\#$ such that
    ${\cal P}'_\# =\{ P'_\# : P_\# \in {\cal P}_\#\}$
     is an $l$-co-connected partition of $M$
    and the correspondence $\mu_\# : {\cal P}_\# \lo {\cal P}'_\#$
    defined by $\mu_\#(P_\#)=P'_\#$ is a matching.
    For every $P \in \cal P$ define
    $P'=\mu(P)=\cup \{ P'_\# : P_\# \subset P \}$
    and ${\cal P}'=\{ P' : P \in {\cal P}\}$.
    Then by 2.2 and 2.3 of \cite{lev} imply that ${\cal P}'$ is
    an  $l$-co-connected partition of $M$ and
    $\mu$ is a matching between ${\cal P}$ and ${\cal P}'$.
    We will call ${\cal P}'$ the modification of ${\cal P}$
    induced by the modification ${\cal P}'_\#$.

    \subsection{  A radial modification}\label{radial}

   Let a cube $B^m=B^q \times B^{m-q}$ be PL-embedded in an $m$-dimensional
    triangulated manifold $M$.
   We consider    $B^m$ with the  linear and coordinate
     structure  induced
    by a  PL-embedding of $B^m$ in
    the Euclidean space $R^m=R^q \times R^{m-q}$ as the unit cube
    $B^m=[-1,1]^m$ with
        $B^q=B^m \cap R^q$ and
    $B^{m-q}= B^m \cap R^{m-q}$.
    For $b \in B^m$, $|b|$ stands for the maximal absolute value
    of the coordinates of $b$.

     Let $r$ be the radial projection
      $r: B^m \setminus O\lo \partial
     B^m$
      and
     $ B_*^m =\frac{1}{2} B^m =\{ |b|\leq 1/2 : b \in B^m \}$.
     The radial modification $P^*$ of a PL-subcomplex $P$ of $M$ is the set

    $P^*=P$ if $P \cap B^m =\emptyset$ and
    $P^*=(P \setminus \Int B^m)\cup
          (r^{-1}(P \cap \partial B^m))$ if
           $P \cap B^m \neq \emptyset$
   \\ and the radial modification
   of a decomposition  $\cal P$ of $M$ is the decomposition
   ${\cal P}^*=\{P^*: P \in {\cal P} \}$ of $M \setminus O$.

   \subsection{ A property of the radial modification}
   \label{property}
   Adopt the notation and the assumptions of \ref{radial}.
   Let  $p: B^m \lo B^q$ be the projection
     to $B^q$. Assume that the decomposition
     ${\cal P}$ of $M$ is such that
     for every $P \in \cal P$, $P\cap B^m =(B^q \cap P) \times B^{m-q}$.

     Suppose that  $f : X \lo Y$ is a map from a dense subspace $X$ of
     $M$
    such that $B^q \subset X$ and  every $x \in B^q$ is
    a regular point of $f$ (see Section \ref{sec-resol}) and
    suppose that ${\cal E}$ is an open cover of $M$.

        Fix  a triangulation
    $\cal T$ of $B^q$ such that $\st^2 \cal T$ refines
    ${\cal E}$.
     Denote  ${\cal B}=\{ \Delta \times B^{m-q} : \Delta \in
     {\cal T} \}$.
     Replacing (if needed)
    $B^{m-q}$ by a smaller cube
     assume that\\

    (*) for every two disjoint
    simplexes $\Delta_1$ and $\Delta_2$ of ${\cal T}$ there
    are disjoint neighborhoods of
    $\Delta_i \times B^{m-q}\subset V_i$ in $M$
    such that
      the closures of $f(V_1 \cap X)$
      and $f(V_2\cap X)$ do not intersect
      in $Y$.\\

  Let us show that there is a PL-homeomorphism
    $h : M \lo M$ such that\\

    (**)  $h$ does not move the points  of
    $(M \setminus \Int   B^m)\cup B^q$, $p(h(c))=p(c)$
    for every  $c \in B^m$ and  $h$ has the property:
     for every $y\in Y$ such that $y$ belongs to the closure of
    $f(h( B^m \setminus \Int B_*^m)\cap X)$  there are a neighborhood
    $G$ of $y$ and $E\in \cal E$
    such that for every $x$  in
     $ h^{-1}( (B^m \setminus\Int B^m_*)\cap f^{-1}(G))$   we have that
    $p(r(x)) \in E$ and
     for every $P \in \cal P$, $P$ intersects $E$ if
     $P^*$ intersects
     $h^{-1}((B^m \setminus \Int B^m_*)\cap f^{-1}(G))$.\\

    Denote ${\cal B}^*=\{ B^* : B \in {\cal B} \}$.
    For  $c=(a, b) \in B^m$   denote
    by $a$ and $b$ the coordinates  of $c$ in $B^q$ and
    $B^{m-q}$ respectively.  Note that if
    $c=(a, b),
    c'=(a', b')\in B^m\setminus O$ are such that $a'=a$ and
    $|b'|=|b|$ then $p(r(c))=p(r(c'))$. Hence
    if  $c=(a, b),
    c'=(a', b')\in B^m\setminus O$ are such that $a'=a$ and
    $|b'|=|b|$ then for every $B^* \in {\cal B}^*$
    we have that $c \in B^*$ if and only if $c' \in B^*$.
    This means that for every  $a \in B^q$ and $B^* \in {\cal B}^*$
    the intersection $B^* \cap (a \times B^{m-q})$ is
    concentric about $a$ (with respect to the norm $|..|$).

    For disjoint $B^*_1, B^*_2 \in {\cal B}^*$
     denote
    $\epsilon(B^*_1, B^*_2)=\inf \{ |c_1 -c_2| : c_i \in B^*_i
    \setminus \Int B_*^m \}$ and let $\epsilon
    =\min\{ \epsilon(B^*_1, B^*_2): B^*_i \in {\cal B}^*,
    B^*_1 \cap B^*_2 =\emptyset\}$. Clearly $\epsilon >0$.
    Let $n$ be a natural number such that $n > 2/\epsilon$.
    Since  $B^q$ consists of   regular points of $f$
    we can find $0< \delta_1< \delta_2 < \dots < \delta_n < 1$
    such that for  $C^+_i = \{ c=(a,b) \in B^m : |b|\geq \delta_i \}$
    and $C^-_i = \{ c=(a,b) \in B^m : |b|\leq \delta_i \}$
    the closure of $f(C^+_i\cap X)$ does not intersect the closure
    of $f(C^-_{i-1}\cap X)$ for every $1< i \leq n$, and
    $\delta_n$ is so close to $0$ that for every $y$ in the closure
    of $f( C_n^- \cap X)$ there is a neighborhood $G$ of $y$
    such that ${\rm diam} p(f^{-1}(G)) \leq \epsilon /2$. Define
    the piece-wise linear function $g : [0,1] \lo [0,1]$
    such that $g(0)=0$, $g(1)=1$ and $g(\delta_i)=i/(n+1), 1 \leq i \leq n$.

    Define a PL-homeomorphism $h : M\lo M$
    such that $h$ does not move the points of $(M \setminus \Int B^m)\cup B^q$,
    for every $c \in B^m$, $p(h(c))=p(c)$ and for every
    $c=(a,b) \in B^m$ such that $|a| \leq n/(n+1)$,
    $h(c)=(a, g(|b|)b)$.
    Denote $C^- =\{ c \in B^m : |c|\leq   n/(n+1)\}$
    and $C^+ =\{ c \in B^m : |c|\geq   n/(n+1)\}$.

    It is easy to see that for
     every disjoint $B^*_1$ and $B^*_2$ in ${\cal B}$,
     the closures
    $f((h((B^*_1\cap C^-)\setminus \Int B^m_*) \cap X)$
    and $f((h((B^*_2\cap C^-)\setminus \Int B_*^m)\cap X)$ do not
    intersect in $Y$. It is also easy to see that since
   $h$ leaves the $B^q$-coordinate of
   the points in $B^m$ unchanged we can
    choose $n$  to be so large and hence $C^+$ to be so close
    to $\partial B^m$ that  (*) will imply that
    for
     every disjoint ${B}_1$ and ${B}_2$ in ${\cal B}$,
     the closures
    $f(h({B^*}_1\cap C^+)\cap X)$ and $f(h({B^*_2}\cap C^+)\cap X)$ do not
    intersect in $Y$.

    Thus   for
     every disjoint $B^*_1$ and $B^*_2$ in ${\cal B}^*$
    we have that the closures
    $f(h((B^*_1\setminus \Int B_*^m)\cap X)$ and
    $f(h(B^*_2\setminus \Int B^m_*)\cap X)$ do not
    intersect in $Y$.

    Take $y $  in the closure of $f(h(B^m \setminus \Int
    B^m_*)\cap X)$.
    The last property  implies that
    there is a neighborhood $G$ of $y$ such that
    $h^{-1}(f^{-1}(G) \cap (B^m \setminus B_*^m))$ is contained in
    an element $C$ of $\st^2{\cal B}^*$. Since
    $r(C)$ is contained in an element  of
    $\st^2 \cal    B$ and $\st^2 \cal T$ refines  $\cal E$
    there is $ E\in \cal E$ such that $p(r(C))) \subset E$.
    Thus we get that if for $P \in {\cal P}$  the radial modification
    $P^*$ of $P$ in ${\cal P}^*$ intersects
      $h^{-1}(f^{-1}(G) \cap (B^m\setminus \Int B^m_*)$
    then $P$ intersects $E$. Clearly  for every
    $x \in h^{-1}(f^{-1}(G)\cap (B^m \setminus \Int B^m_*))$,
     $p(r(x)) \in E$ and (**) follows.

     \subsection{Combining  radial and black hole modifications}\label{comb}

     Let $M$ be  an $m$-dimensional
     triangulated manifold,
     $F$ a PL-subcomplex of $M$ and
     ${\cal P}$  a decomposition of $M$ which forms a partition on
     $M \setminus F$. Assume that
     $B^m =B^q\times     B^{m-q} \subset M$ is a ball constructed
     for  creating or improving connectivity of a finite intersection of
     ${\cal     P}$, see 2.5 and 4.2 of \cite{lev} and
     \ref{create} of this section.

     Assume that $B^m$ is represented as in \ref{radial}.
     For a PL-subcomlex $A$ of $M$ define
     $A' = (A^* \setminus \Int B_*^m)\cup A_*$
     where $A^*$ is the radial modification of $A$ and
     $A_* = \frac{1}{2}(A \cap B^m)=
     \{ \frac{1}{2}a : a \in A \cap B^m    \}$.
     It is easy to see that ${\cal P}'=\{P': P \in {\cal P}\}$
     is a decomposition of $M$ which forms  a partition
     on $M\setminus F'$ and we preserve in ${\cal P}'$
     restricted to $M \setminus \Int B^m_*$ the connectivity
     properties
     that we have in $\cal P$ restricted to $M \setminus \Int B^m$.
     Notice that ${\cal P}^*$ and ${\cal P}'$ both
      restricted to  $M \setminus \Int B^m_*$ coincide and
       ${\cal P}'$ restricted
     to $B^m_*$ is an exact copy of $\cal P$ restricted to
     $B^m$.
    Then the block hole modification used in see 2.5 of \cite{lev}
     for improving connectivity can be carried out in $B^m_*$
     instead of $B^m$ leaving ${\cal P}'$ on $M \setminus \Int B^m_*$
     as the radial modification.
      Thus we can shift the block hole modification
      to $B^m_*$ combining it with the radial
      modification.

    \end{section}

    \begin{section}{ Proof of Theorem \ref{shrink}}\label{sec-shrink}

     Assume that $X$ is homeomorphic to $M(k)$ of an
    $m$-dimensional $(n-1)$-connected triangulated manifold $M$ with
    $m\geq 2n+1$ and
     $k=m-n-1$. Let ${\cal P}$ be
     a decomposition
    of an open subset ${ M}_\#$ of $M$ containing
    $X$.
    We say that $y \in Y$ is
    a $\cal P$-regular point
     if
     there is neighborhood $G$ of $y$ in $Y$ such that
    $G \cap X $ is contained in an element of
    $\st{\cal P}$,
    and
   we say that  $y\in Y$ is $\cal P$-singular otherwise.
   Note that the set of $\cal P$-singular points is
   closed in $Y$ and does not meet $X$ and therefore it is
   a $Z$-set in $Y$.

    Assume that  $Z' \subset Z  \subset  Y$,
    $Z'$ and $Z$ are homeomorphic to $n$-dimensional
    \no spaces, $Z$ is a $Z$-set in $Y$, $Z'$ is closed
    in $Z$, a map
    $g_Z: Z\lo Z'$ is a $UV^{n-1}$-retraction and
    a map $g_Y : Y \lo Z$ is a retraction. Consider
    an open cover ${\cal W}_{Z'}$ of $Z'$ and denote
    by ${\cal W}_M$ the extension of
    $g_Y^{-1}(g^{-1}_Z({\cal W}_{Z'})|X$ to $M$. Replacing
    $M$ by an open subset of $M$ containing $X$ we assume that
    ${\cal W}_M$ covers $M$.

    Note that each time when we replace
    $M$ by
     a smaller open set containing $X$ we automatically replace
     everything that was already defined on $M$ (sets, covers,
     decompositions
     etc.) by its restriction to that smaller open set.

     By  a nice partition we mean
   a rational $(m-n+1)$-co-connected partition with no
   non-empty intersections of $\dim \leq m-n-1$. Note that
   every non-empty
   finite intersection of a nice partition has
   a dense subset lying in $X$.

    By  a neighborhood  we always mean an open neighborhood.
     If $Q$ is  an open subset
    in one of the spaces $X, Y$ or $M$ then  the closure $\cl Q$
     of $Q$ is considered
     with respect to that space where $Q$ is open.

   \begin{proposition}
   \label{shift}
    {\rm {\bf (Shifting Singularities)}}
   Assume that  $X$ and $Y$ as in Theorem \ref{shrink}
   and  $M$,  $Z'$,  $Z$, $g_Z$, $g_Y$,
   ${\cal W}_{Z'}$ and ${\cal W}_M$
   as above.

   Then for every
   open neighborhood $H$ of $Z$ in $Y$  there is
   a smaller neighborhood $Z \subset Q_Z  \subset H$
   of $Z$ in $Y$
   such that for every nice
    partition ${\cal P}$ of an open subset $M_\#$ of
    $M$ such that $X \subset M_\#$
   and  the set of ${\cal P}$-singular
   points of $Y$
    is contained in $Q_Z$ we can do the following:

   for every neighborhood
   $Q_{Z'}$ of $Z'$ in $Y$
   and a ball $K $ of
   $\dim \leq n$ PL-embedded in $M$  such that
   $K \subset X$, the manifold $M_\#$ can be modified
   to an open subset $M'_\#$ of $M_\#$ and each $P \in {\cal P}$
   to a PL-subcomplex $P'$ of $M'_\#$ such that

   (1) $X \subset M'_\#$, ${\cal P}'=\{ P' : P \in {\cal P}\}$
   is a nice partition of $M'_\#$,
   the correspondence $P \lo P'$ is a matching
   between ${\cal P}$ and ${\cal P}'$;

   (2)
    $P \cap((X\setminus H) \cup K)=
   P'  \cap ((X\setminus H) \cup K)$ for every
   $P \in \cal   P$   (that is
   ${\cal P}$ and ${\cal P}'$ coincide on
   $(X\setminus H) \cup K$);

   (3)
   the set of ${\cal P}'$-singular points of $Y$
   is contained in $Q_{Z'}$ and

    (4) for every $P \in {\cal P}$,
    $P'$ is contained in $\st(P, {\cal W}_M)$.

   \end{proposition}
   Let us show how Theorem \ref{shrink} can be derived from
   Proposition \ref{shift}.
   \\\\
   {\bf Proof of Theorem \ref{shrink}.}
   Consider the space
   constructed in Proposition \ref{p2.ad} and  denote
   $Z_i =\cup_{j \geq i } X_j$. Let
   $h_i : A_i \lo A_0$ be a homeomorphism such that
   $h_0$ is the identity map.  By
   Proposition \ref{p2.3}, extend $h_i$ and $h_{i+1}$ to
   a $UV^{n-1}$-map $\psi_i : X_i \lo A_0$ and
   define the $UV^{n-1}$-retraction
   $g_{i}^Z: Z_i \lo Z_{i+1}$ by $g_i^Z(x)=h_{i+1}^{-1}(\psi_i(x))$
   for $x \in X_i$ and the map $g^A_i : Z_i \lo A_0$ by
    $g^A_i(x)=\psi_j(x)$ for $x \in X_j, j\geq i$. Note
    that $g^A_i=g^A_{i+1}\circ g_i^Z$.

     Let ${\cal V}'_Y$ and ${\cal V}'_X$  be open
     covers of $Y$ and $X$ respectively  such that
   $\st^2 {\cal V}'_Y$ refines ${\cal V}_Y$
   and  $\st^2 {\cal V}'_X$ refines ${\cal V}_X$.
   Replacing $M$ by
   an open subset of $M$ containing $X$ we assume that
     the extension  ${\cal V}^M_Y$ of
   ${\cal V}'_Y|X$   to $M$ and
   the extension  ${\cal V}^M_X$ of
   ${\cal V}'_X$ to $M$ both
   cover $M$.

   Take a nice partition
   ${\cal P}$ of $M$ such that $\st^2 \cal P$ refines both
   ${\cal   V}^M_Y$ and
   ${\cal   V}^M_X$.
    Denote by $T({\cal P })$ the set of ${\cal P}$-singular points of
   $Y$ and recall that $T({\cal P })$ is a $Z$-set in $Y$.
    Embed  $A_0$ in $Y$ as a $Z$-set in $Y$
    such that
    $T({\cal P}) \subset A_0$.
    Let  ${\cal V}''_Y$ be an open cover of $Y$
    such that $\st^4 {\cal V}''_Y$  refines ${\cal V}'_Y$.
     Consider a sequence
    ${\cal W}^A_i$ of open covers of  $A_0$ such that
    $\st^2{ \cal W}^A_{i+1} $ refines $ {\cal W}^A_{i} $
    and ${\cal W}^A_{0} $ refines ${\cal V}''_Y$.
     Via the map
    $g^A_0 : Z_0 \lo A_0$ extend the embedding of
    $A_0$ to a $Z$-embedding  $Z_0 \subset Y$ such that
    $(g^A_0)^{-1}({\cal W}^A_0)$ refines $\st{\cal V}''_Y$.
    Let $g_Y : Y \lo Z_0$ be any retraction and
    let
    a neighborhood  $H$ of $Z_0$ be
    such that
    ${\cal W}^Y_0=g^{-1}_Y((g^A_0)^{-1}({\cal W}^A_0))$ restricted to $H$
    refines $\st^2{\cal V}''_Y$. Then
    $\st{\cal W}^Y_0$ restricted to $H$
    refines ${\cal V}'_Y$.

    Choose a neighborhood $H_i$ of $Z_i$ in $Y$ such that
     $H_0=H$,
    the closure of $H_{i+1}$ is contained in $H_i$ and
    $\cap_{i=0}^{\infty} H_i =\emptyset$.
    Denote $H_i^M$= the extension of $H_i$ to $M$.
     Replacing $M$
    by an open subset of $M$ containing $X$ we assume
    that the closure $\cl{H}_{i+1}^M$ of  $H_{i+1}^M$
    is contained in
     $H^M_i$ and $\cap_{i=0}^{\infty} H_i^M =\emptyset$.
     Denote
    ${\cal P}_0=\cal P$, $M_0=M$,
     ${\cal W}_0^{M}=$ the extension of
     ${\cal W}^Y_0$
     to $M$ and replacing $M$ by an open subset of $M$ containing
     $X$ assume that ${\cal W}_0^{M}$ covers $M_0$. For every $i>0$
     fix an open subset $M_i$ of $M$ such that $X \subset M_i$,
      $M_{i+1} \subset M_i$
     and the extension
    ${\cal W}_i^{M}$ of ${\cal W}^Y_i=g_Y^{-1}((g^A_0)^{-1}({\cal W}^A_i))$ to $M_i$
    covers
    $M_i$. Denote $g^Y_0=g_Y$, $g^Y_i=g^Z_{i-1}\circ \dots\circ g^Z_0 \circ g_Y : Y \lo Z_i$
    for $i >0$
    and ${\cal W}^{Z}_i=(g^A_{i+1})^{-1}({\cal W}^A_i)$. Note
    that
    ${\cal W}_i^Z$ is an open cover of $Z_{i+1}$ and
     ${\cal W}^Y_{i}=(g^Y_i)^{-1}((g^Z_i)^{-1}({\cal
    W}^Z_{i}))$.

    By Proposition \ref{shift} there is a neighborhood $Q_i$
    of $Z_i$ such that
    the conclusions
     of  Proposition \ref{shift} are satisfied
     for $M$,$H$, $Z$, $Z'$, $g_Y$,
     $g_Z$, ${\cal W}_{Z'}$ and ${\cal W}_M$
 replaced by $M_i$, $H_i$, $Z_i$, $Z_{i+1}$,$g^Y_i$, $g^Z_i$,
 ${\cal W}^Z_{i}$ and
  ${\cal W}_{i}^{M}$
  respectively. Note that $Q_i \subset H_i$.

   Let $K_i \subset X$ be
     a sequence of
    balls
     PL-embedded in $M$ which will be chosen later.
     Recall that $T({\cal P}_0) \subset Z_0$ and
     therefore $T({\cal P}_0) \subset Q_0$.
  Set $M'_0=M_0$ and  apply   Proposition \ref{shift} to
     construct for every $i$
  a nice partition ${\cal P}_i $ of an open subset $M'_i$ of $M$
   such that
  $X \subset M'_i \subset M_i$, the set $T({\cal P}_i)$ of
  ${\cal P}_i$-singular points of $Y$
  is contained in $Q_i$, and ${\cal P}_{i+1}$ and $M'_{i+1}$
  are produced by Proposition \ref{shift} as the output
  ${\cal P}'$ and $M'_\#$ when the input
  $\cal P$, $M_\#$,  $Q_{Z'}$ and $K$ corresponds to ${\cal P}_i$,
   $M'_i$,  $Q_{i+1}$ and $K_i$ respectively with the other parameters
   as above. Note that replacing
   $M'_\#$ by $M'_\#\cap M_{i+1}$ we indeed can assume
   that $M'_{i+1} \subset M_{i+1}$.
   Let $\mu^i : {\cal P}_i \lo {\cal P}_{i+1}$
   be the matching sending each element to its modification
   and let $\mu_i : {\cal P}={\cal P}_0 \lo {\cal P}_i$
   be the composition of the corresponding chain of matchings.
   Denote $P_i =\mu_i(P)$ for $P \in {\cal P}$. Since $\mu_i$
   is a matching, $\mu_i(P) $ is well defined for every
   finite intersection $P$ of ${\cal P}$.

   Define $M'=\cup_{i=0}^\infty (M'_i \setminus \cl   H^M_{i})$,
    $P'= \cup_{i=0}^\infty (P_i \setminus \cl   H^M_{i})$ for
   $P \in {\cal P}$.
    It is clear that $M'$ is an open
   subset of $M$ containing $X$, $Y$ has no
    ${\cal P}'$-singular point
   and ${\cal P}'=\{ P' : P \in {\cal P}\}$ is
   a rational partition of $M'$. However we can say nothing about
   the connectivity of ${\cal P}'$ and it may well happen
   that  $P'=\emptyset $  for $P \in {\cal P}$. In order to cope
   with this problem we need to choose the balls $K_i$
   in the following manner.
   Take a bijection $\alpha : \N \lo \N \times \N \times \N $
   such that for  $\alpha(s)=(j, t, l) $ we have
   $s \geq \max \{j, t, l\}$.
   Arrange the non-empty finite intersections of
   ${\cal P}$ into a sequence $P^0,P^1,P^2, \dots$.   For every
   $j$ and $t$ such that $\mu_{j} (P^{t})
   \setminus \cl H^M_{j} \neq \emptyset$
   choose a countable collection
   ${\cal F}(j,t)=\{f(j,t,l): l=0,1,2,\dots\}$ of maps
   from spheres of
   $\dim \leq \dim P^{t} -m + n-1$ into $(\mu_{j} (P^{t}) \setminus
   \cl H^M_{j})\cap X$ such that ${\cal F}(j,t)$
    generates the homotopy groups of
   $\mu_{j} (P^{t}) \setminus \cl H^M_{j}$ in $\dim \leq
   \dim P^{t} -m +   n-1$. For
   $j$ and $t$ such that $\mu_{j} (P^{t}) \setminus
    \cl H^M_{j} = \emptyset$
   define ${\cal F}(j,t)=\emptyset$.
   Set  $i(0)=0$ and
   for each $s$ define  by induction a set $K'_{s} \subset X$  and
   $i(s)\geq s$ as follows.
      Let $\alpha(s)=(j,t,l)$. If
      ${\cal F}(j,t)=\emptyset$ define
      $K'_{s}$ to be any singleton in $\mu_{i(s)} (P^{t})\cap X$.
      If ${\cal F}(j,t)\neq\emptyset$ then
      by \ref{general} let
      $K'_{s}$  be a ball of $\dim \leq n$ PL-embedded in
      $\mu_{i(s)}(P^{t})$ such that $K'_{s} \subset X$
      and the map $f(j, t, l)$ is null-homotopic
      in $(\mu_{i(s)}(P^{t})\setminus \cl H^M_{i(s)}) \cup K'_{s}$.
      Recall that $i(s)\geq s \geq j$ and therefore the image
      of $f(j,t,l)$ is contained in $\mu_{i(s)}(P^{t})\setminus \cl
      H^M_{i(s)}$.
      Since
    $\cap_{i=0}^\infty  \cl H^M_i =\emptyset$
    there is $i(s+1)> i(s)$ such that $K'_{s} \cap \cl H^M_{i(s+1)} =\emptyset$.
     Now define
     $K_{i}= K'_{s}$ for $i(s) \leq i < i(s+1)$.
     Note  that   $K_{i}$ is already defined  after ${\cal P}_i$
     is constructed and
     therefore  the construction of ${\cal P}_{i+1}$
     is well-defined.
    Then, by (2) of Proposition \ref{shift},
    we get  that $P'$ will be non-empty for
    every non-empty  $P \in {\cal P}$, ${\cal P}'$ is a nice partition
    of $M'$ and  the correspondence
    $\mu : {\cal P} \lo {\cal P}'$  defined
    by $\mu(P)=P'$ is a matching of partitions.

    Let $P \in \cal P$. Recall that
    $P=P_0$. By (4) of Proposition \ref{shift} we have that
     $P_{i+1} \subset \st(P_i, {\cal W}_{i}^{M})$ and since
     $\st^2 {\cal W}_{i+1}^{M}$ refines ${\cal W}_i^{M}$
     we have that $P' \subset \st(P_0,  \st{\cal W}_0^{M})=
     \st(P,  \st{\cal W}_0^{M})$. Since $P \setminus  H^M_0=
     P_0 \setminus H^M_0=P'\setminus H^M_0 $  we get
     $P' \subset (P\setminus H^M_0)\cup
     ( \st(P_0,  \st{\cal W}_0^{M})\cap H^M_0) $ and
     since
      $\st{\cal W}^Y_0$ restricted to $H=H_0$
      refines ${\cal V}'_Y$
       we have
      that $P' \subset \st (P, {\cal V}^M_Y)$.
      Recall that ${\cal P}$ refines ${\cal V}^M_Y$
      and $\st {\cal V}'_Y$ refines ${\cal V}_Y$.
      Hence  for every $P \in \cal P$ there is
      an element of
       ${\cal V}_Y$ which contains both $P \cap X$
       and $P'\cap X$.

       Using
      Theorem  \ref{unknot} construct
      for every finite intersection $P$ of $\cal P$
       (by induction on $\dim P$)
       a homeomorphism
      $h_P : P\cap X \lo P'\cap X$ such that
      $h_P$ agree on common intersections.
      Then the  homeomorphisms $h_P$ define
      the corresponding  homeomorphism
      $h : X \lo X$ such that for every
      $P \in \cal P$, $h(P\cap X) =P'\cap X$ and hence
       $x$ and $h(x)$ are  ${\cal V}_Y$-close
      for every $x \in X$.

      Recall that $\st \cal P $ restricted  to $X$
       refines ${\cal V}_X$ and $Y$ has
       no ${\cal P}'$-singular point. Then for every
       $y$ in $Y$ there is a neighborhood $G$
       such that $G\cap X$ is contained in an element
       of $\st {\cal P}'$ and hence
       $h^{-1}(G\cap X)$ is contained in an element
       of $\st \cal P$. Therefore
       $h^{-1}(G\cap X)$ is contained in an element
       of ${\cal V}_X$ and  the theorem follows.
       \hfill $\Box$\\
       \\
       For proving Proposition \ref{shift} we need the following
       \begin{lemma}\label{lemma}
       Let $X$ and $Y$ be as in Theorem \ref{shrink}.
       Then there are a $Z$-embedding $g : X\lo Y$
       of $X$ into $Y$ and a continuous retraction
       $f : Y \lo g(X)$ such that both $g$
       and $f$ can be chosen to be arbitrarily close
       to the inclusion of $X$ into $Y$ and
       the identity map of $Y$ respectively.
       \end{lemma}
        {\bf Proof.}
        Let us first show that
         there are a $Z$-embedding $\phi : X \lo X$
   and a continuous retraction $r : X \lo \phi(X)$ such
   that both $\phi$ and $r$ can be chosen to
   be  arbitrarily
   close to the identity map of $X$. Indeed, take
   any topological copy $X'$ of $X$ such that
   $X'$ is $Z$-embedded in $X$ and let
   $\psi_1 :X_1= X \lo X'$ be any homeomorphism.
    By Proposition \ref{p2.3}
   there is a $UV^{n-1}$-retraction
   $\psi_2 : X_2=X \lo X'$.
   Take any $Z$-embedding $\phi : X_1 \lo X_1$
   which is sufficiently close to the identity map.
    Then, by Theorem \ref{unknot1}, the homeomoprhism
    $\psi_1 \circ \phi^{-1} : A_1= \phi(X_1) \lo A_2=X' \subset X_2$
     extends
   to a homeomorphism $h : X_1 \lo X_2$ such that
   the retraction
   $r= h^{-1} \circ \psi_2 \circ h
    : X_1 \lo \phi(X_1) $  is close to the identity map of $X_1$.

    Note that the inclusion of $X$ into $Y$ is
    a $UV^{n-1}$-map.
     Choose  a $Z$-embedding $\xi :\phi(X) \lo Y$
    to be so close to the inclusion of $\phi(X)$ into $Y$
    that, by Proposition \ref{p2.1}, the map
    $\xi^{-1} :\xi (\phi(X)) \lo \phi(X)$ extends
    to a map $\eta : Y \lo X$ which is
     close to the identity map of
    $Y$. Then $f: Y \lo \xi(\phi(X))$ defined by
    $f(y)=\xi (r (\eta(y)))$ is close to the identity
    map of $Y$ provided $r$ is close  to
    the identity map of $X$. Denote $g=\xi \circ \phi$
    and we are done.
     \hfill $\Box$\\\\
      {\bf Proof of Proposition \ref{shift}.}
      Fix a (sufficiently large)
       natural number $\omega$ which will be determined
      later and which  depends only on $n$.
       We assume that indices involving $i$ and $j$
      run from $1$ to $\omega$. Thus, for example, if we
      write $i+1$ or $j-2$ as an index we assume that $1 \leq i+1 \leq \omega$
      and $1 \leq j-2 \leq \omega$. All the properties with
      indices involving   $i$ and $j$
      are automatically restricted
     to the cases when the indices make sense (=remain within
     the range from $1$ to $\omega$).

      Recall that  $Z'$ and $Z$ are  \no spaces and $g_Z$ is $UV^{n-1}$.
      Then there are
            open  covers ${\cal W}_{Z'}(i)$ of
        $Z'$ and  neighborhoods $Q_Z(i)$ of $Z$ in $Y$  such that
         $\st^2{\cal W}_{Z'}(i) \prec_{n-1} {\cal W}_{Z'}(i+1) $,
         ${\cal W}_{Z'}(\omega)\prec {\cal W}_{Z'}$,
        $\cl Q_Z(i)\subset Q_Z(i+1)$, $Q_Z(i) \subset H$  and
         for the open
        covers
        ${\cal W}_Y(i)= g^{-1}_Y(g^{-1}_Z({\cal W}_{Z'}(i)))$
        of $Y$
        we have that
       $\st {\cal W}_Y(i)$ restricted to
        $Q_Z(i)$ is an $(n-1)$-refinement  of ${\cal W}_Y(i+1)$
        restricted $Q_Z(i+1)$.

        The sets $Q_Z(i)$
        can be constructed starting from
        $Q_{Z}(\omega)$  and
        choosing for each $i= \omega-1,\dots, 1$
        the set $Q_{Z}(i)$ to be so close to $Z$ that
        every map $\alpha$ from a sphere $S$ of $\dim \leq n-1$
        to $Q_{Z}(i)$ can be ${\cal W}_Y(i)$-closely
        homotoped   inside $Q_{Z}(i+1)$ to the  map
        $g_Y \circ \alpha : S \lo Z$.

        We are going to show that
        $Q_Z =Q_Z(1)$ satisfies the conclusions of the proposition.

         Let   $M_\#$ be an open subset of $M$ containing $X$,
       ${\cal P}$ a nice partition of $M_\#$,
        $Q_{Z'}$ a neighborhood of $Z'$ in $Y$
        and  $K$  a ball  of
   $\dim \leq n$ PL-embedded in $M$  such that
   $K \subset X$.
        Denote by
        $T({\cal P})$ the set of the ${\cal P}$-singular points
        of $Y$
        and assume that  $T({\cal P})$ is contained in $Q_Z$.
        Since our choice of $Q_Z$ does not depend
        on $M$ we can
        replace $M$ by $M_\#$ and assume that
        $\cal P$ is a partition of $M$. Note
        that without loss of generality we  can always replace $M$
         by any open subset of $M$ containing $X$
         (with the automatic replacement of $\cal P$ by
        the restriction of $\cal P$ to that subset).

        Let $Q_{Z'}$ be a neighborhood of $Z'$ in $Y$.
        Since $Z'$  is a $Z$-set in $Y$ there are
                    neighborhoods
        $Q_{Z'}(i)$ of $Z'$ in $Y$  so
         that
              $ Q_{Z'}(i) \subset Q_{Z'} $, $ Q_{Z'}(i) \subset
              Q_Z(1)$,
       $\cl Q_{Z'}(i)\subset Q_{Z'}(i+1)$,
        ${\cal W}_Y (i)$ restricted to $Q_{Z'}(i)$
        is an $(n-1)$-refinement of
         ${\cal W}_Y (i+2)$ restricted to $Q_{Z'}(i+2)$
         and for every $i > j $ we have that
       ${\cal W}_Y (i)$ restricted to
        $Q_{Z'}(i) \setminus \cl Q_{Z'}(j)$ is
        an $(n-1)$-refinement of ${\cal W}_Y (i+2)$
        restricted to  $Q_{Z'}(i+2) \setminus \cl Q_{Z'}(j-2)$
        and  for every point in $Z'$ there is
        a  ${\cal W}_Y(i)$-close point contained in
        $Q_{Z'}(i) \setminus \cl Q_{Z'}(j)$.

        The sets $Q_{Z'}(i)$  can be constructed
        starting from  $Q_{Z'}(\omega)$.
        Assume that $Q_{Z'}(\omega), \dots Q_{Z'}(i+1)$ are
        already constructed. Construct $Q_{Z'}(i)$ as follows.
         The first property of $Q_{Z'}(i)$ can be taken care of
        in a way similar to the one used
        for  constructing $Q_{Z}(i)$. In order
        to satisfy the second property
         choose
        a $Z$-embedding $\beta : Y \lo Y$ such that
        $\beta(Y) \cap Z' =\emptyset$ and
        $\beta$ is so close to the identity map of
        $Y$ that there is a neighborhood
        $Q_{Z'}(i)$ of $Z'$ such that
        $\cl Q_{Z'}(i) \subset Q_{Z'}(i+1) \setminus \beta(Y)$
        and
         every map  from a sphere
        of $\dim \leq n-1$ to
         $Y \setminus \cl Q_{Z'}(i+1)$
         can be sufficiently closely homotoped
         in $Y \setminus \cl Q_{Z'}(i)$
         to  a map to $\beta(Y)$. If $\beta$
         is  close enough  to the
         identity map of $Y$ then $Q_{Z'}(i)$
         has the required properties.

        Denote by $Q_{Z}^M(i)$, $Q_{Z'}^M(i)$ and ${\cal W}_Y^M(i)$
        the extensions to $M$ of
        $Q_{Z}(i)$, $Q_{Z'}(i)$ and ${\cal W}_Y(i)$ all restricted to
        $X$.
        Replacing $M$ by an open subset of $M$ containing $X$
         assume that ${\cal W}_Y^M(i)$ cover $M$,
        $\cl Q_{Z}^M(i)\subset Q_{Z}^M(i+1)$ and
        $\cl Q_{Z'}^M(i)\subset Q_{Z'}^M(i+1)$.

       Let
         $A$ be a space  of $\dim \leq n$ and
         $B$ a closed subset of $A$. Consider maps
       $\alpha: A \lo Y$  and $\beta : B \lo Y$ and
       let  $\delta> 5n+5$.
       From
        the properties of $Q_{Z'}(i)$, $Q_Z(i)$ and ${\cal W}_Y(i)$
        one can derive the following.
        If the image of $\beta$ is contained in
        $Q_Z(i)$ and in an element of ${\cal W}_Y(i)$ then
        $\beta$ can be extended to $\beta' :L \lo Y$
        such that the image of $\beta'$ is contained in
        $Q_Z(i+\delta)$ and in an element of ${\cal
        W}_Y(i+\delta)$. If the image of $\beta$ is contained
        in $ Q_{Z'}(i+3\delta) \setminus \cl Q_{Z'}(i+2\delta)$
        then $\beta$ extends to a map into
        $ Q_{Z'}(i+4\delta) \setminus \cl Q_{Z'}(i+\delta)$.
         If the image of $\beta$ is contained
        in $Q_{Z'}(i+3\delta) \setminus\cl Q_{Z'}(i+2\delta)$
        and $\alpha$ is any extension of $\beta$ then there
        is an extension
         $\alpha' : L \lo Q_{Z'}(i+4\delta)
         \setminus \cl Q_{Z'}(i+\delta)$ of $\beta$ such that
         $\alpha$  and $\alpha'$ are ${\cal W}_Y(i+4\delta)$-close.
        \\\\
        {\bf Constructing an initial subdivision ${\cal P}_\#$
        and an initial  modification ${\cal P}_\#(1)$. }\\
     Since $T({\cal P}) \subset Q_Z(1)$ and $T({\cal P}) \cap K
     =\emptyset$ there are open neighborhoods $Q_T(i)$ of
     $T({\cal P})$ in $Y$ such that $\cl Q_T(i) \subset
     Q_Z(1)\setminus K$ and $\cl Q_T(i) \subset Q_T(i+1)$. Then
     there are  open covers ${\cal W}_T(i)$ of $Y$  such that
     $\st^2 {\cal W}_T(i) \prec_{n-1} {\cal W}_T(i+1)$,
     $\st(\cl Q_T(j),{\cal W}_T(i)) \cap \st (Y\setminus
     Q_T(j+1), {\cal W}_T(i) )=\emptyset$,
      $\st(\cl Q_Z(j), {\cal W}_T(i)) \cap
      \st(Y\setminus Q_Z(j+1), {\cal W}_T(i))=\emptyset$,
     $ {\cal W}_T(\omega)\prec {\cal W}_Y(1)$,
     and $\st^2 {\cal W}_T(\omega)$ restricted to
     $X\setminus Q_T(1)$
     refines $\st{\cal P}$.

     Replacing $M$ by an open subset of $M$ containing $X$
     assume
     that the extension to $M$ of ${\cal W}_T(1) | X$ covers
     $M$. By \ref{subdiv} there are $(m-n+1)$-co-connected
     rational partitions  ${\cal D}(i)$ of $M$ such that
     ${\cal D}(i)$ is a subdivision of ${\cal D}(i+1)$,
     $\st^3 {\cal D}(i)\prec_{n-1}\st{\cal D}(i+1)$,
     ${\cal D}(\omega)$ is a subdivision of ${\cal P}$ and
     ${\cal D}(\omega)$ refines the extension to $M$ of
     ${\cal W}_T(1) | X$.
     Deleting from $M$ the finite
     intersections  of ${\cal D}(i)$ of $\dim \leq m-n-1$
     we assume that ${\cal D}(i)$ are nice partitions of $M$.

     By   Lemma \ref{lemma} there is a $Z$-embedding $g : X \lo Y$
     and a continuous retraction $f : Y \lo g(X)$ such that
     $g$ and $f$  are ${\cal W}_T(1)$-close to the inclusion
     of $X$ into $Y$ and the identity map of $Y$ respectively.
     Denote $F_Y(D) = f^{-1}(g(D\cap X))$ for $D \in {\cal D}(1)$
     and ${\cal F}_Y = \{ F_Y(D) :  D \in {\cal D}(1) \}$.
     Note that ${\cal F}_Y$ is a closed locally finite cover
     of $Y$. Replacing $M$ by an open subset of $M$ containing
     $X$ we assume that for
     ${\cal F}_M = \{ F_M(D) :  D \in  {\cal D}(1) \}$,
     where
     $F_M(D)=$the closure of $F_Y(D) \cap X$
     in $M$,
     we have that ${\cal F}_M$
     is a locally-finite cover  of $M$  such that
     $F_M(D_1)\cap F_M(D_2) =\emptyset$ provided
      $D_1 \cap D_2=\emptyset, D_1, D_2 \in {\cal D}(1)$.

     Let ${\cal P}_\#$ be a rational $(m-n+1)$-co-connected partition
     of $M$ such that ${\cal P}_\#$ is a subdivision
     of ${\cal D}(1)$ and
     $\st(F_M(D_1), {\cal P}_\#)\cap
     \st( F_M(D_2),{\cal P}_\#) =\emptyset$ provided
      $D_1 \cap D_2=\emptyset, D_1, D_2 \in {\cal D}(1)$.
      Deleting from $M$ the finite intersections of
      ${\cal P}_\#$ of $\dim \leq m-n-1$ we assume that
      ${\cal P}_\#$ is a nice partition of $M$.

     For every $P \in {\cal P}_\#$ that intersects
     $X \setminus  Q_T(\omega)$ define $P(1)=P$.
     Arrange into a sequence $P_1, P_2,\dots$ the elements
     of ${\cal P}_\#$ that do not intersect $X \setminus  Q_T(\omega)$
     and define by induction on $s$ the set
     $P_s(1)$ as follows: $P_1(1)=$the union of the elements
      of ${\cal P}_\#$ that do not intersect $X \setminus  Q_T(\omega)$
      and do intersect $F_M(P_1)$, and
      $P_{s+1}(1)=$the union of the elements
      of ${\cal P}_\#$ that do not intersect
       $(P_1(1)  \cup \dots \cup P_s(1))\cup (X \setminus  Q_T(\omega))$
      and do intersect $F_M(P_{s+1})$.

      Define ${\cal P}_\#(1)=\{ P(1) : P \in {\cal P}_\#\}$.
      Clearly ${\cal P}_\#(1)$ covers $M$. However we may have
      that $P_s(1)=\emptyset $ for some $P_s $.
      To avoid such situation  modify the elements
      of ${\cal P}_\#(1)$ as follows. For every $P_s$  such that
      $P_s(1)=\emptyset $ choose a point
      $x_s$ in $X$ such that $x_s$ is contained in $\Int P_{x_s}$
      of some $P_{x_s} \in {\cal P}_\#$,
       $x_s$ is
      ${\cal P}_\#$-close to some point of $F_M(P_s)$
      and the set  of all the points $x_s$ is discrete in $X$.
      Replacing $M$ by a smaller open subset containing
      $X$ we assume that the set of $x_s$ is discrete
      in $M$ as well. Let
      $B_{x_s}$ be an  $m$-dimensional
      PL-ball rationally embedded in $\Int P_{x_s}$
      such that $ x_s \in \Int B_{x_s}$,
       $\Int P_{x_s} \setminus B_{x_s}\neq \emptyset$
       and the collection of the balls $B_{x_s}$
       is discrete in $M$.
       Then cutting the interior of the balls $B_{x_s}$ from  the
       elements of ${\cal P}_\#(1)$ and defining
       $P_s(1)= B_{x_s}$ we get that
       ${\cal P}_\#(1)$ is a rational partition
       of $M$ for which $\mu_1: {\cal P}_\# \lo {\cal P}_\#(1)$,
       defined by $\mu_1(P)=P(1)$, is a $0$-matching.

       Define ${\cal A}(i)$ as the cover of $M$
       consisting of
       the elements of $\st{\cal D}(i)$
       and the sets $\st({W}\cap X, {\cal D}(i))$ for
       ${W }\in {\cal W}_T(i)$ such that
       ${W }\cap (Y\setminus  Q_T(\omega -i)) \neq \emptyset$.
       From the properties of  ${\cal W}_T (i)$
       and ${\cal D}(i)$ it follows that\\

        ($C_1$)
       \\
              each element of  ${\cal A}(i)$ is
              a union of elements of ${\cal P}_\#$,
       $ {\cal A}(i)$ refines $\st{\cal P}$ for $i \leq \omega-3$,
      ${\cal A}(i)$ refines
      ${\cal W}^M_Y(i+2)$,
       $\st {\cal A}(i) \prec_{n-1} {\cal A}(i+2)$
       and
       $\st(\cl Q_Z^M(i), {\cal A}(i)) \subset Q_Z^M(i+2)$.\\\\
          {\bf Constructing   modifications ${\cal P}_\#(i)$.}
          We are going to construct modifications
          ${\cal P}_\#(i)$ of ${\cal P}_\#$ and   rational subcomplexes
          $F(i)$ of $M$
           such  that ${\cal P}_\#(i)$
          is a rational decomposition of $M$ forming a partition
          on $M \setminus F(i)$  and
          ${\cal P}_\#$ admits a one-to-one correspondence
          $\mu_i :{\cal P}_\# \lo {\cal P}_\#(i)$ to
          ${\cal P}_\#(i)$ where $P(i)=\mu_i(P)$ is the modification
          of
          $P \in {\cal  P}_\#$ in ${\cal P}_\#(i)$.

         Note that during the construction
          we will often replace $M$ by  smaller open subsets
          containing $X$. Our goal is to get in the end of the
          construction a nice partition ${\cal P}(i)$ for which
            $\mu_i$ is a matching.
            We will carry out the construction such that for
            every   $i$
             there will exist $j \geq i$ for which
            the following  important  condition holds.\\

            ($C_2$)
            \\
             $F(i) \cap ((X\setminus Q_Z(j))\cup K)=\emptyset$,
            ${\cal P}_\#(i)$ and ${\cal P}_\#$ coincide on
            a neighborhood in $X$ of
            $(X\setminus Q_Z(j))\cup K$ (that is
            the intersections of $P$ and $P(i)$ with that
            neighborhood coincide for every $P \in {\cal P}_\#$),
             $P(i)=P$ if
            $P\in {\cal P}_\#$ is contained in $  M
            \setminus Q_Z^M(j)$,
            for every $P \in {\cal P}_\#$ there is an element
            of ${\cal W}_Y(j)$ containing both $P\cap X$ and $P(i)\cap X$,
            and for every $y \in Y \setminus
            Q_{Z'}(j)$ there is a neighborhood $G$ of $y$ in $Y$
            such that $\mu_i^{-1}(\st (G\cap X , {\cal P}_\#(i))$
            is contained in an element of ${\cal A}(j)$. \\
           \\

       Let us verify  that the initial modification
        ${\cal  P}_\#(1)$  satisfies ($C_2$)
       with $F(1)=\emptyset $ and $j=4$.
       The condition concerning $K$ is satisfied because
         $(X\setminus Q_Z(4))\cup K \subset X \setminus
       \cl Q_T(\omega) $. It is also clear that
       $P(1)=P$ if $P \subset M \setminus Q_Z^M(4)$. Since
       $f$ and $g$ are
       ${\cal W}_T(1)$-close to the inclusion of $X$ into $Y$
       and the identity map of  $Y$ respectively  we have
        for
       $P \in {\cal P}_\#$  that $P \cap X$ and $P(1)\cap
       X$ are contained in an element of $\st^4 {\cal W}_T(1)$
       and,
       since  ${\cal W}_T(1)$ refines ${\cal W}_Y(1)$,
       $P \cap X$ and $P(1)\cap
       X$ are contained in an element of ${\cal W}_Y(4)$.

       Let $y \in Y$ be such that
        $f(y) \in Y \setminus Q_T(\omega -1)$. Then for every element $G$
       of ${\cal W}_T(1)$ containing $y$ we have that
       $\mu_1^{-1}(\st (G\cap X , {\cal P}_\#(1))$ is contained
       in an element of $\st^5 {\cal W}_T(1)$ which
       intersects $Y\setminus Q_T(\omega-4)$ and
       therefore it is contained in an element of
        ${\cal A}(4)$.

       Now let  $y \in Y$ be such that
       $f(y) \in  Q_T (\omega-1)$.
       Take a neighborhood of $y$ of the form
       $G =$the interior in $Y$ of $g^{-1}(V)$
       where $V =\st(f(y), g({\cal D}(1)|X))$.
       Then no element of ${\cal P}_\#$ intersecting
       $Q_T(\omega)$ will intersect $G$ and hence
              for $P \in {\cal P}_\#$, $P(1)$ intersects
       $G$ only if $P \subset \st(g^{-1}(f(y)), \st^2{\cal D}(1))$.
       Thus  $\mu_1^{-1}(\st (G\cap X , {\cal P}_\#(1))$ is contained
       in an element of $\st^5 {\cal D}(1)$ and therefore it
       is contained in ${\cal A}(4)$. Condition $(C_2)$
       has been verified for ${\cal P}_\#(1)$.\\

       The purpose  of the modifications ${\cal P}_\#(i)$ is
       to gradually improve the level of connectivity
       and the level of matching to ${\cal P}_\#$. It will be done
       following the construction \ref{level}. This construction
       is a combination of 2.7 and 2.9 of \cite{lev} with
       a use of \ref{create} for creating intersections.

       Note that we do not need to
        worry   about the rationality  of ${\cal P}_\#(i)$
        and $F(i)$  inside the constructions themselves.
        Indeed, assume that ${\cal P}_\#(i)$
        and $F(i)$ are constructed such that $(C_2)$ holds.
        By \ref{general}
         we may assume
         that every finite intersection of ${\cal P}_\#(i)$
         has a dense subset in $X$.
         Now replacing $M$ again
         by an open subset containing $X$
         and  we can choose by \ref{rem-mov}
         a PL-homeomorphism   $h: M \lo M$ to be sufficiently
         close to the identity map of $M$
         such that $h$ does not move the points on
         the extension to $M$ of  a neighborhood of
         $(X \setminus Q_Z(j))\cup K$ in $X$ on which
         ${\cal P}_\#(i)$ and ${\cal P}_\#$ coincide,
         $h$ moves  ${\cal P}_\#(i)$ to
         a rational decomposition and $F(i)$ to a rational
         subcomplex of $M$ such that
          ${\cal P}_\#(i)$ and $F(i)$ replaced by
         $h({\cal P}_\#(i))$  and $h(F(i))$ respectively
          will satisfy $(C_2)$ with $j$ replaced by $j+1$.

         Let us analyze separately the constructions involved
       in \ref{level}.
        We  always  assume that the finite intersections of
       ${\cal P}_\#(i)$ have dense subsets on $X$.
         Fix   $0 \leq s \leq n-1$, $i$ and $j$,
         set $l=m-s+1$, $q=s+1$
         and
       assume that $(C_2)$ and the following
       condition hold for $i$, $j$,  $s$ and $l$.\\

       $(C_3)$\\
      $\dim F(i) \leq l-2$,
      ${\cal P}_\#(i)$ is $l$-co-connected
       on $M \setminus F(i)$, $\mu_i$ induces an $s$-matching
       between ${\cal P}_\#$ and ${\cal P}_\#(i)$ restricted to
       $M \setminus  F(i)$. \\
      \\

      Let ${\cal G}$ be a cover  of $Y \setminus  Q_{Z'}(j)$ by
      open sets in $Y$ such that for every $G \in G$,
      $\mu_i^{-1}(\st (G\cap X , {\cal P}_\#(i))$
            is contained in an element of ${\cal A}(j)$.
            Let  ${\cal E}_Y$ be
            the collection of open subsets of $Y$
            consisting of the elements of ${\cal G}$
            and the set $Q_{Z'}(j)$.
            Denote by ${\cal E}$ the extension to $M$ of
            ${\cal E}_Y$ restricted to $X$. Replacing $M$ by an open
            subset containing $X$ assume that $\cal E$
            covers $M$. Choosing $\cal G$ to be sufficiently fine
            assume that
            $\st(\cl Q_{Z'}^M (j), \st^6{\cal E})
           \subset Q_{Z'}^M(j+1)$.
\\

             {\bf 1). Improving connectivity of
       intersections.}  Fix $0\leq t \leq m-l+1$ and
        assume that the finite intersection of
        ${\cal P}_\#(i)$ restricted to
        $ M\setminus F(i)$ of $\dim > m-t$ are
        $(l-1)$-co-connected.
       We are going to improve the
        connectivity of the intersections of $\dim =m-t$
         using the construction 2.7 of \cite{lev}.

        Adopt the notation of 2.7 of \cite{lev} with
        ${\cal P}$ and $F$ replaced by ${\cal P}_\#(i)$ and
        $F(i)$ respectively.
           Recall that in   2.7 of \cite{lev} we choose
           for every finite intersection $P$ of
           ${\cal P}_\#(i)$ such that $\dim P \setminus F(i)=m-t$
           and $P\setminus F(i)$ is not $(l-1)$-co-connected, a countable
           collection  ${\cal F}_P$ of PL-embeddings
           of an $(q-t-1)$-dimensional sphere into
           $P \setminus F(i)$ such that ${\cal F}_P$
           represents all the elements of the
           $(q-t-1)$-dimensional homotopy group
           of $P \setminus F(i)$. After that we extend all the embeddings
           in ${\cal F}_P$ for all $P$
           to a collection ${\cal F}^{\partial \Delta}$
            of  PL-embeddings
             $f_{\partial \Delta} : \partial B^q \lo \Int M$
        of  a $q$-dimensional PL-ball $B^q$ such that
        the images of the embeddings in ${\cal F}^{\partial \Delta}$
        form a discrete family in $M \setminus R^{\partial \Delta}$
        for
          a PL-presented subset $R^{\partial \Delta}$ of $M$ with
         $\dim R^{\partial \Delta} \leq q-1$
        (here we may consider $\Delta$ and $\partial \Delta$
        only as  a  part of the notation).
         By \ref{rem-discr} we can assume that $R^{\partial \Delta}$
         is rationally presented
         and replace $M$ by $M \setminus R^{\partial \Delta}$.
        By \ref{improv}
        we  assume that
        $f_{\partial \Delta} (\partial B^q) \subset X$ for
        every $f_{\partial \Delta} \in {\cal F}^{\partial
        \Delta}$.

        Recall that
         $\st(\cl Q_Z^M(j), {\cal P}_\#) \cap
         \st(M \setminus Q_Z^M(j+1), {\cal P}_\#)=\emptyset $.  By
         $(C_2)$, the elements of  ${\cal
         P}_\#$ that are contained in
         $M\setminus   Q_Z^M(j)$
         are left unchanged in ${\cal P}_\#(i)$. Then, since ${\cal P}_\#$ is
         a nice partition, the finite intersections
         of ${\cal P}_\#(i)$ that are not $(l-1)$-co-connected
          intersect $ Q_Z^M(j)$. The embeddings in
         ${\cal F}^{\partial \Delta}$ needed for improving
         the connectivity of
         a finite intersection $P$ of ${\cal P}_\#(i)$ lie
         in $\st(P, {\cal P}_\#(i))$. Thus we may assume that
         the images
         of the embeddings in ${\cal F}^{\partial \Delta}$
         refine $\st{\cal P}_\#(i)$ restricted to $Q_Z^M(j+1)$
         and therefore they will refine $\st {\cal W}_Y(j+1)$
         restricted to $Q_Z(j+1)$.

         Consider $B^q$ as   the unit cube $B^q=[-1,1]^q$
         in the $q$-dimensional
         Euclidean space and let $r : B^q\setminus O \lo \partial
         B^q$ be the radial projection.
         Extend each $f_{\partial
         \Delta}$ to a map
         $f^{+}_{  \partial \Delta} : B^q \setminus
         \Int \frac{3}{4} B^q  \lo
          X$ such that
          $f^{+}_{\partial \Delta} \circ r$ and
          $f_{\partial  \Delta}$ are both ${\cal W}_Y(j)$
          and ${\cal E}_Y$-close,
      the images of  the maps $f^+_{\partial \Delta}$
      are contained in $Q_Z(j+3)$ and
       the images of  the maps $f^+_{\partial \Delta}$
      restricted to $\partial \frac{3}{4} B^q$
     form a discrete family in $Y$ (recall that $Y \setminus X$
     is a $\sigma$-$Z$-set in $Y$).
     Note that
     the images of  the maps $f^+_{\partial \Delta}$
     refine ${\cal W}_Y(j+3)$.
\\\\
        Set $\delta = 10n+10$ and
        let $r' : \frac{3}{4} B^q \setminus O \lo
        \partial \frac{3}{4} B^q $ be the radial projection.
         Denote

        $B_0 (f^+_{\partial \Delta})=
        {r'}^{-1}
        [(f^+_{\partial \Delta})^{-1}[\cl Q_{Z'}(j+3\delta)
        \setminus Q_{Z'}(j+2\delta)] \cap
        \partial \frac {3}{4}B^q] \cup \frac{1}{2} B^q$,

        $B_- (f^+_{\partial \Delta})=
        {r'}^{-1}
        [({f}^{+}_{\partial \Delta})^{-1}[\cl Q_{Z'}(j+3\delta)
        ] \cap
        \partial \frac {3}{4}B^q] \setminus \Int \frac{1}{2} B^q $,

        $B_+ (f^+_{\partial \Delta})=
        {r'}^{-1}
        [({f}^{+}_{\partial \Delta})^{-1}[Y
        \setminus Q_{Z'}(j+2\delta)] \cap
        \partial \frac {3}{4}B^q] \setminus \Int \frac{1}{2} B^q$.
\\\\
      By the properties of $Q_{Z}(j), Q_{Z'}(j)$ and
     ${\cal W}_Y(j)$  the maps  $f^+_{ \partial\Delta}$
     can be extended
     to  maps $f'_{\Delta} : B^q \lo X$
     such that  the images  of $f'_{\Delta}$ are contained in
     $Q_Z(j+5\delta)$ and refine $W_Y(j+5\delta)$,
     the images  of $f'_{\Delta}$ restricted to
     $ \frac{3}{4} B^q$ form
     a discrete family in $Y$ and for every
     $f'_{\Delta}$ we have that
      the image  of $f'_{\Delta}$ restricted to
     $B_0 (f^+_{\partial \Delta})$ is contained in
     $ Q_{Z'}(j +4\delta) \setminus \cl Q_{Z'}(j+\delta)$,
      the image  of $f'_{\Delta}$ restricted to
     $B_+ (f^+_{\partial \Delta})$ is contained in
     $ Y \setminus \cl Q_{Z'}(j+\delta)$
     and  the image  of $f'_{\Delta}$ restricted to
     $B_- (f^+_{\partial \Delta})$ are contained in
     $ Q_{Z'}(j +4\delta) $.

    By \ref{general} we can replace $M$ by an open subset
    containing $X$ and assume that
    the images of  the maps $f'_{\partial \Delta}$
      restricted to $\partial \frac{3}{4} B^q$
     form a discrete family in $M$.
    By 2.6 of \cite{lev},
     the maps $f'_\Delta$ can be arbitrarily
     closely approximated by
     maps $f''_\Delta : B^q \lo \Int M  \setminus R^\Delta$
     such $f''_\Delta$ coincides with $f'_\Delta$
     on $\partial B^q$ and
     the images  of $f''_\Delta$  are contained
     and form a discrete family in
     $M  \setminus R^\Delta$ for some
     PL-presented subset
     $ R^\Delta$ of $M$
     with
      $\dim R^\Delta \leq n$. By \ref{rem-discr}
     assume that $R^\Delta$ is rationally presented
     and replace $M $ by $M \setminus R^\Delta$.
     Then by \ref{rem-dig} each  map $f''_\Delta$
     can be arbitrarily closely approximated
     by a PL-embedding $f_\Delta : B^q \lo \Int M $
     such that $f_\Delta$ coincides with $f_{\partial \Delta}$
     on $\partial B^q$, $f_\Delta$
      extends to a PL-embedding
       $f_B$   of $ B^m =B^q \times B^{m-q} $
      into $M $ with the properties
      described in 4.2 of \cite{lev}. Then if
      $f_\Delta$ is  sufficiently
      close to $f''_\Delta$ we can choose $f_B$ so
      that the images of $f_B$
      form a discrete family in
      $M $ and the sets
      $f_B((\frac{3}{4}B^q) \times B^{m-q}) \cap X$
      form a discrete  family in $Y$. By the reasoning
      similar to the one of \ref{improv} there is
       a point $a \in \Int B^{m-q}$
       arbitrarily close to
      $O  \in B^{m-q}$  such that
       $f_B(B^q \times a )\subset X\setminus K$. Replacing
      $B^{m-q}$ by a small cube centered at $a$ we assume that
        $f_B(B^q ) \subset X\setminus K$, and
      $f_B$ and $f_B \circ p$ are $\cal E$-close where
      $p : B^m \lo B^q$ is the projection.

      Thus without loss of generality we may assume  that
       the images of $f_\Delta$ are contained in $X$ and
       all the above   properties
       hold when $f'_\Delta$, $f^+_{\partial \Delta}$
       and $f_{\partial \Delta}$ are
           replaced by $f_\Delta$,
         $f_\Delta$
       restricted to $B^q \setminus \Int \frac{3}{4}B^q$
       and $f_\Delta$ restricted to $\partial B^q$
       respectively
       with the only change that because the images
        of $f_\Delta$ may slightly move
       and
       the distance between $f^{+}_{\partial \Delta} \circ r$ and
          $f_{\partial  \Delta}$  may slightly increase
          on the set $B^q \setminus \Int \frac{3}{4}B^q$
          where we do not have the discreteness
          in $Y$ of the images
          we will require that
          the images of $f_\Delta$ are contained in
          $Q_Z(j+6\delta)$ and refine ${\cal W}_Y(j+6\delta)$,
          and
            $f^{+}_{\partial \Delta} \circ r$ and
          $f_{\partial  \Delta}$ are $\st{\cal W}_Y(j)$ and
           $\st{\cal E}_Y$-close.

    In addition,  for each $f_B$ we can replace $B^{m-q}$ by a
    smaller cube and  assume that
   $f_B(B_+(f^+_{\partial \Delta}) \times B^{m-q})\subset M
    \setminus \cl Q^M_{Z'}(j+\delta)$ and
     $f_B(B_+(f^-_{\partial \Delta}) \times B^{m-q})\subset
    Q^M_{Z'}(j+4\delta)$, the images of $f_B$ are contained
    in $Q_Z^M(j+6\delta)$ and refine ${\cal W}_Y^M(j+6\delta)$,
    and the images of $f_B$ are outside
     some neighborhood  of $K$ in $M$.

    For a given embedding $f_B$ identify $B^m=B^q\times B^{m-q}$ with
    $f_B(B^m)$. Clearly
    replacing $B^{m-q}$ by a smaller cube we may
    assume that (*) of \ref{property}
    is satisfied with $f$ being the inclusion of
    $X$ into $Y$.
     Now apply   the radial
    modification to ${\cal P}_\#(i)$ as described in \ref{radial},
    and let   $h: M \lo M$ be a PL-homeomorphism
    satisfying (**) of \ref{property} with $\cal P$
    replaced by ${\cal P}_\#(i)$ and $f$ replaced by the inclusion
    of $X$ into $Y$. Note that now $r$ stands for
    the radial projection $r : B^m \setminus O\lo \partial B^m$
    which extends the  radial  projection for $B^q$ previously denoted
    by $r$.

    Since the images of $f_B$ are discrete in $M$
    this procedure can be done independently for
    each $f_B$. Let us denote by ${\cal P}_\#^*(i)$
    the result of all the radial modifications
    (note that we also modify $F(i)$) and
    let us denote again by $h$ the resulting
     PL-homeomorphism  of $M$ for all embeddings $f_B$.

     Let $y \in Y\setminus \cl Q_{Z'}(j+4\delta)$. Let us
     show
     that\\

     $(C_4)$
\\
     there are a neighborhood $G$ of $y$ in $Y$
     and an element $E$ in $\cal E$ lying
     in $M \setminus \cl     Q_{Z'}^M (j+1)$
     such that for every $P \in {\cal P}_\#(i)$,
     $ h(P^*)\cap X$ intersects $G$ only
     if $P $ intersects $\st(E, \st^5{\cal E})$
     where
      $P^*\in {\cal P}_\#^*(i)$ is the modification
      of $P$.\\

      Denote $F_Y=$ the closure in $Y$ of the union of
     $f_B((\frac{3}{4} B^q ) \times B^{m-q})\cap X$
     for all $f_B$ and recall that
     the family of the sets
     $f_B((\frac{3}{4} B^q ) \times B^{m-q})\cap X$ is discrete
     in $Y$.

      Consider first the case
     $y \notin F_Y$. Take a neighborhood $G$ of $y$ in $Y$ such
     that $G \cap F_Y =\emptyset$ and $G$ is contained in
     an element $E_Y$ of ${\cal E}_Y$ lying in $Y\setminus \cl Q_{Z'}(j+1)$.
     Let $E\in \cal E$ be the extension of
     $E_Y \cap X$ to $M$ and let $P \in {\cal P}_\#(i)$ be such
     that
     $h(P^*)\cap X$ intersects $G$.
     If $P\cap X$ intersects $G$ then clearly
     $P$ intersects $E$.
     Assume that $P\cap X$ does not intersect $G$.
       Take
      $x \in G\cap( h(P^*) \cap X)$. Then
      there is $f_B$ such that $x$ belongs  to $f_B(B^m)$.
      Identify $f_B(B^m)$ with $B^m$ and note that
     $r(h^{-1}(x)) \in P$ and,
      $x$ and $r(h^{-1}(x))$ are  $\st^3 {\cal E}$-close
      because $G \cap F_Y =\emptyset$.
      Hence $P$ intersects $\st(E, \st^3{\cal E})$ and
      $(C_4)$ follows.

      Now let $y \in F_Y$
      and $y \in Y\setminus \cl Q_{Z'}(j+4\delta)$. Then there is
      (only one) $f_B$ such that
      $y $ belongs to
      the closure in $Y$ of $ f_B((\frac{3}{4} B^q ) \times B^{m-q})\cap
      X$. Identify  $B^m$ with $f_B(B^m)$.
      Take a neighborhood $G$ of $y$ in $Y$
       such that $G\subset Y \setminus \cl Q_{Z'}(j+4\delta)$
       and $G$ satisfies (**) of \ref{property}
      for  an element $E \in \cal E$.  Clearly $G$ can be replaced
      by any smaller neighborhood of $y$.

      Assume   $y$ does not belong
      to the closure in $Y$ of $\partial B^m \cap X$. Replace $G$ by
      a smaller neighborhood and assume that $G \cap X \subset B^m$.
      Let
       $x \in (G \cap X)\cap ((\frac{3}{4} B^q ) \times B^{m-q})$,
        $x'=p(r((h^{-1}(x)))$, $x''=x'$ if $x'\in
        \frac{3}{4}B^q$ and
         $x''=\frac{3x'}{4|x'|}$ otherwise
        (note that $O$ does not belong to $G \cap X$
        since $ (\frac{1}{2})B^q \times B^{m-q}
        \subset Q_{Z'}^M(j+4\delta)$) .
       From the properties of $f_B$ it follows that
       $x'' \in B_+(f^+_{\partial \Delta})
       \subset
       M \setminus \cl  Q_{Z'}^M (j+\delta)$.
      By (**) of \ref{property} we have that $x' \in E$.
      Then,  since $x' $ and $x''$ are
      $\st^2\cal E$-close, we have
       that $E$ lies in
      $M \setminus \cl Q_{Z'}^M(j+1)$ and $(C_4)$ follows.

       Now assume that
      $y$ belongs to the closure  of $\partial B^m \cap X$ in $Y$.
      Replace $G$ by a smaller neighborhood of $y$ and assume
      that $G$ does not intersect
      $F_Y \setminus (B^m \cap  X)$ and is contained in an element
      of ${\cal E}_Y$. Take a point $x \in G \cap B^m$ which is
      so close to $\partial B^m$ that $p(x)$ and $p(r(h^{-1}(x))$
      are $\cal E$-close. Then $G\cap X$ is contained in
      $\st(E, \st{\cal E})$. Let $P \in {\cal P}_\#(i)$ be such that
      $P^*\cap (G\cap X)\neq \emptyset$.
       By (**) of \ref{property}
      we have that if $h(P^*)$ intersects $B^m \cap (G \cap X)$ then
      $P$  intersects $E$. If $P^*$ does not intersect $B^m$ then
      by the reasoning similar  to the one applied in the case
      $y \notin  F_y$ we conclude that $P$ intersects
      $\st (E, \st^5{\cal E})$. Since $G\cap X$ is contained in
      $\st(E, \st{\cal E})$ and $G\subset Y \setminus Q_{Z'}(j+\delta)$
      we have that $E \subset M \setminus \cl Q_{Z'}(j+1)$
      and $(C_4)$ follows.
      Thus $(C_4)$ has been  verified in
      all cases.

       Now we are ready to
      construct ${\cal P}_\# (i+1)$ and $F(i+1)$.
      Fix $f_B$ and identify $B^m$ with $f_B(B^m)$.
      Modify ${\cal P}_\#(i)$ on $B^m$ as described
      in \ref{comb} and combine this modification with
      the block hole modification  shifted to
      $B^m_*=\frac{1}{2} B^m$. Denote by ${\cal P}^{**}_\#(i)$
      and  $F^{**}(i)$ the result of all these modifications
      for all $f_B$.
      Denote $\Omega=$the union
      of $f_B( (\frac{1}{2} B^q )\times B^{m-q})$ for all $f_B$.
      Note that ${\cal P}^*_\#(i)$ and
      ${\cal P}^{**}_\#(i)$ coincide on $M \setminus
      \Omega$ and
      $\Omega$ is contained in $Q_{Z'}^M (j +4\delta)$.
      Define ${\cal P}_\# (i+1)=h({\cal P}^{**}_\#(i))$
      and $F(i+1)=h(F^{**}(i))$. Note that in $(C_4)$,
      $\st (E, \st^5{\cal E}) \subset M\setminus  Q_{Z'}^M (j)$
      since $E \subset M \setminus \cl Q_{Z'}^M(j+1)$ and
      $\st (\cl Q_{Z'}^M (j), \st^6{\cal E}) \subset
      Q_{Z'}^M(j+1)$. Then $(C_4)$ implies that for every
      $y \in Y \setminus \cl Q_{Z'}(j + 4\delta)$ there is
      a neighborhood $G$ of $y$ in $Y$ such that
      $\mu^{-1}_{i+1}(\st (G\cap X , {\cal P}_\#(i+1)))$
      refines $\st^6 {\cal A}(j)$.

      Since the images of $f_B$ are contained in
      $Q_{Z}^M(j+6\delta)$,  do not intersect $K$
      and form a discrete family in $M$ we have
      that $F(i+1) \subset Q_{Z}^M(j+6\delta)$
      and there is
      a neighborhood  of
      $( M \setminus Q_{Z}^M(j+6\delta))\cup K$ in $M$
      on which ${\cal P}_\#(i)$ and ${\cal P}_\#(i+1)$
      coincide.
     Thus we finally  get that ${\cal P}_\#(i+1)$ and
      $F(i+1)$ satisfy $(C_2)$ with $j$ replaced
      by $j +6\delta$, ${\cal P}_\#(i+1)$ and
      $F(i+1)$
        satisfy $(C_3)$   and
      all the intersections of $\dim \geq m-t$  of
      ${\cal P}_\#(i+1)$ restricted to $M \setminus F(i+1)$
      are $(l-1)$-co-connected.\\

         {\bf 2). Creating intersections.} Assume that ${\cal
          P}_\#(i)$ is $(l-1)$-co-connected on $M \setminus F(i)$
          and satisfies
          $(C_2)$ and $(C_3)$.
          Following \ref{level}   we need to create
          the missing intersections of ${\cal P}_\#(i)$.
          Let $P=P_0 \cap\dots \cap P_{s+1}$ be a non-empty
          intersection of distinct elements of ${\cal P}_\#$
          such that $(P_0(i) \cap \dots \cap P_{s+1}(i)) \cap (M
          \setminus F(i))=\emptyset$. Since ${\cal P}_\#$
          and ${\cal P}_\#(i)$ coincide on $M \setminus Q_Z^M(j)$
          we have that $P_0(i), \dots, P_{s+1}(i) \subset
          Q_Z^M(j+1)$.  Since for every $P\in {\cal P}_\#$
           there is an element of   ${\cal W}_Y^M(j)$
          containing both $P$ and $P(i)$, the union of
           $P_0(i),\dots,P_{s+1}(i)$ is contained
           in an element of ${\cal W}_Y^M(j+2)$.
           Then, by \ref{create}, the maps ${\cal F}^{\partial
           \Delta}$ can be chosen so that the images of
            ${\cal F}^{\partial   \Delta}$ are
          contained in $Q_Z^M(j+2)$ and
             refine
          ${\cal W}_Y^M(j+2)$. The rest of the
          construction is identical to the previous construction
          of improving connectivity of intersections.\\

       {\bf 3). Absorbing simplexes.} This is the last
       construction in \ref{level}.  Assume that
        ${\cal  P}_\#(i)$ is $(l-1)$-co-connected on $M\setminus F(i)$,
       ${\cal  P}_\#(i)$
        satisfies
          $(C_2)$ and $(C_3)$ and all the intersections of
          ${\cal  P}_\#(i)$ (restricted to $ M\setminus F(i)$)
          are brought from ${\cal P}_\#$ (by $\mu_i$).
          Following \ref{level}  add  to $F(i)$ all the finite
          intersections of ${\cal P}_\#(i)$ which are not
          brought from the finite intersections of ${\cal P}_\#$
          of $\dim \leq m-s-1$. Since ${\cal P}_\#$ and
          ${\cal P}_\#(i)$
          coincide on $M \setminus Q_Z^M(j)$ we have
          $F(i) \subset Q_Z^M(j+1)$.
          Since
           ${\cal P}_\#$ and
          ${\cal P}_\#(i)$
          coincide on a  neighborhood of $K$ in $M$
          (because they coincide on a  neighborhood of $K$ in $X$)
          we have that $F(i) \cap K =\emptyset$.

          Following 2.9 of \cite{lev} fix a sufficiently
            fine triangulation of $M$ underlying ${\cal P}_\#(i)$
          and $F(i)$  We will first absorb the $(m-s-1)$-simplexes
          of $F(i)$ lying in $M \setminus Q_{Z'}^M(j)$.
          Replacing $M$ by an open subset containing $X$
          and the triangulation of $M$ by a finer triangulation
          assume that for
           a simplex $\Delta$ of $\dim =m-s-1$  lying
           in $F(i)\setminus Q_{Z'}^M(j)$ we have that
           $\mu_i^{-1}(\st(\Delta, {\cal P}_\#(i)))$
           is contained in an element $A$ of ${\cal A}(j)$.
           Take an element $A'$ of ${\cal A}(j+1)$ such that
           $A \subset A'$ and the inclusion of $A $ into $A'$
           induces the zero-homomorphism of the homotopy
           groups in $\dim \leq n-1$. Recall that ${\cal P}_\#$
           underlies both $A$ and $A'$. Denote
           $F=$ the union of the finite intersections of
           ${\cal P}_\#$ of $\dim < m-s-1$ which are not
           brought from ${\cal P}_\#(i)$ restricted to
           $M \setminus F(i)$.
           Then $\mu_i$ induces a matching between
           ${\cal P}_\#$ restricted to $M \setminus F$
           and ${\cal P}_\#(i)$ restricted to $M \setminus F(i)$.
            Since $\dim F <
           m-s-1$ we have that ${\cal P}_\#$
           is $(l-1)$-co-connected  the inclusion
           $A\setminus F \subset A'\setminus F$ induces
           the zero-homomorphism of the homotopy groups
           in co-dimensions$\geq l-1$.

            Hence,  by \ref{rem-absorb},
            the construction for absorbing  all simplexes of
            $\dim=m-s-1$ lying in $M\setminus Q_{Z'}^M(j)$
             can be carried out  such that the modification
              ${\cal P}_\#(i+1)$ of ${\cal P}_\#(i)$ will
              refine
              ${\cal B}=\mu_i(\st^2{\cal A}(j+1))$ and
              for every
              $P(i) \in {\cal P}_\#(i)$,
              $P(i+1) \subset \st (P(i), {\cal B})$
              where $P(i+1)\in{\cal P}_\#(i+1)$ is the modification of
              $P(i)$.
             Note that  we also modify $M$, however by \ref{rem-discr}
             the PL-subcomlex of $\dim \leq n-1$ needed to be removed
             from $M$ can be rational and hence $M$ can be
             replaced
             by an open subset containing $X$.
              Also note that by \ref{rem-absorb} we can assume
              that
               ${\cal P}_\#(i+1)$ coincides
               with ${\cal P}_\#(i)$ on a small neighborhood
               of $K$ in $M$.
               Clearly ${\cal P}_\#(i+1)$ coincides
              with ${\cal P}_\#(i)$ outside $\st(Q_Z^M(j+1),
                 {\cal   B})$.
                    Hence ${\cal P}_\#(i+1)$ refines
              ${\cal W}_Y^M(j+10)$ and coincides with ${\cal
              P}_\#$ on $M \setminus Q_Z^M(j+10)$.

               Let $y \in M \setminus Q_{Z'}(j)$ and let $G$ be
   neighborhood of $y$ in $Y$ such that that
   $\mu^{-1}_i(\st(G \cap X, {\cal P}_\#(i))$ is contained in an element
    of
   ${\cal A}(j)$. Then $\mu^{-1}_{i+1}(\st(G \cap X, {\cal P}_\#(i+1))
   \subset \mu^{-1}_{i}(\st(G \cap X, {\cal B})
   \subset \st(\mu^{-1}_i[\st(G\cap X, {\cal P}_\#(i))],
   \mu^{-1}_i{\cal B})\subset\st(\mu^{-1}_i[\st(G\cap X,{\cal
   P}_\#(i))],
   \st^2{\cal A}(j+1))$ and hence
   $\mu^{-1}_{i+1}(\st(G \cap X, {\cal P}_\#(i+1))$ is
   contained in an element of ${\cal A}(j+10)$.

   Thus replacing $i+1$ by $i$ and $j+10$ by $j$
   we can assume that
   ${\cal P}_\#(i)$ and $F(i)$ satisfy $(C_2)$
   and  $\dim F(i) \setminus Q_{Z'}^M(j)< m-s-1$.
   Now will absorb all the simplexes of $F(i)$ of
   $\dim=m-s-1$. Similarly to the previous case we
   define ${\cal B}$. Note that ${\cal B}$ refines
   $W_Y^M(j+10)$ and the elements of ${\cal  B}$ are contained
   in $Q_Z^M(j+10)$.
   Let $\alpha :S\lo M$ be a map from a sphere $S$ of dimension $s$
   such that $\alpha (S)$  lies in $Q_{Z'}^M(j)\setminus F(i)$
   and in an element $B$ of ${\cal B}$, and
   $\alpha$ is null-homotopic in $B \setminus F(i)$.
   Set $\delta =10n+10$.
   Then, since $\dim F(i) \setminus Q_{Z'}^M(j)< m-s-1$,
   it follows from
     the properties of $Q_{Z'}^M(i)$ that  $\alpha$  will be
   null-homotopic in $[[\st (B , {\cal W}_Y^M(j +2\delta))] \cap
   Q_{Z'}^M(j+2\delta)]\setminus F(i) $.

   Thus the construction of absorbing the $(m-s-1)$-dimensional
   simplexes of $F(i)$ can be carried out inside
   $Q_{Z'}^M(j+3\delta )$ such that the modifications
   ${\cal P}(i+1)$ and $F(i+1)$ will satisfy $(C_2)$
   with $j$ replaced by $j+3\delta$ and $(C_2)$ with
   $l$ replaced by $l-1$ and $s$ replaced by $s+1$
    (we take care of $K$ and
    the modification of $M$ similarly to the previous case).
   \\
   \\
    {\bf Constructing the final modification ${\cal P}'$.}
    The construction of the modifications of ${\cal P}_\#$ ends up
    with a rational decomposition
    ${\cal P}_\#(i)$ and a rational subcomplex $F(i)$ of $M$
    such that  ${\cal P}_\#(i)$ and  $F(i)$ satisfy $(C_2)$
    for some $j$,
     ${\cal P}_\#(i)$ forms an $(m-n+1)$-co-connected
     partition on $M \setminus F(i)$, $\dim F(i) \leq m-n-1$
     and $\mu_i $ is an $n$-matching between $ {\cal P}_\# $
     and ${\cal P}_\#(i)$ restricted to $M \setminus F(i)$ .
     Removing  $F(i)$ and the finite intersections of ${\cal P}_\#(i)$
     of $\dim \leq m-n-1$
     from  $M$ we can assume that ${\cal P}_\#(i)$ is a nice
     partition of $M$ and $\mu_i$ is
     a matching. Note that the number of the modifications
     of ${\cal P}_\#$ needed to get ${\cal P}_\#(i)$ depends only
     on $n$ and the increase in the value of $j$ for each
     modification is bounded by a number depending only on $n$.
     Thus we can estimate the maximal value of $j$ and
     assign this value to $\omega$ in the beginning of the proof.

     Following \ref{subdiv} for each $P \in {\cal P}$ define
     $P'=\cup \{ \mu_i(P_\#) : P_\# \in {\cal P}_\# , P_\# \subset P\}$
     and  let ${\cal P}' =\{ P' : P \in {\cal P}\}$.
     Since $\cal P$, ${\cal P}_\#$ and
     ${\cal P}_\#(i)$  are nice partitions, ${\cal P}_\#$
     is a subdivision of ${\cal P}$ and $\mu_i$ is a matching
     between ${\cal P}_\#$ and
      ${\cal P}_\#(i)$ we conclude that ${\cal P}$ is a nice
      partition of $M$ and $\mu : {\cal P} \lo {\cal P}'$
      defined by $\mu(P) =P'$ is a matching. From $(C_2)$
      it follows that $P' \subset \st(P, {\cal W}_Y^M(j))$
      and, since ${\cal W}_Y^M(j)$ refines ${\cal W}_M$
      we get  $P' \subset \st(P, {\cal W}_M)$.
       It also follows from $(C_2)$ that ${\cal P}$ and
       ${\cal P}'$ coincide on a neighborhood of
       $(M \setminus Q_Z^M(j))\cup K$ and since
       $Q_Z(j) \subset H$
       the requirement (2) of the proposition is satisfied.

       Let us verify the requirement (3). Take a point
       $y \in Y \setminus Q_{Z'}$. Then $y \in Y \setminus
       Q_{Z'}(j)$. By $(C_2)$ there is a neighborhood
       $G$ of $y$ in $Y$ such that
       $\mu_i^{-1}(\st (G \cap X, {\cal  P}_\#(i)))$ is contained
       in an element of ${\cal A}(j)$
       and hence, by $(C_1)$,
       $\mu_i^{-1}(\st (G \cap X, {\cal  P}_\#(i)))$ is contained in
       an element of $\st{\cal P}$. Then, since
       $\mu_i$ is a matching and ${\cal P}_\#$ is
       a subdivision of ${\cal P}$ we have
       that $\mu_i(\mu_i^{-1}(\st (G \cap X, {\cal  P}_\#(i))))=
       \st (G \cap X, {\cal  P}_\#(i))$ is contained
       in an element of $\st{\cal P}'$ and therefore $y$
       is not a ${\cal P}'$-singular point. The requirement (3)
       has been verified.
        \hfill $\Box$\\\\

    \end{section}

Department of Mathematics\\
Ben Gurion University of the Negev\\
P.O.B. 653\\
Be'er Sheva 84105, ISRAEL  \\
e-mail: mlevine@math.bgu.ac.il\\\\
\end{document}